\documentstyle[12pt]{article}
\input amssym.def
\begin{document}
\def \Z{\Bbb Z}
\def \C{\Bbb C}
\def \R{\Bbb R}
\def \Q{\Bbb Q}
\def \N{\Bbb N}
\def \P{\Bbb P}
\def \bZ{{\bf Z}}
\def \TZ{{1\over T}{\bf Z}}
\def \wt{{\rm wt}\;}
\def \fg{\frak g}
\def \mod{{\rm mod}\;}
\def \pf{{\bf Proof}.}
\def \x{{\bf x}}
\def \y{{\bf y}}

\def \wt{{\rm wt}}
\def \span{{\rm span}}
\def \Res{{\rm Res}}
\def \End{{\rm End}}
\def \Hom{{\rm Hom}}
\def \<{\langle} 
\def \>{\rangle}
\def \be{\begin{equation}\label}
\def \ee{\end{equation}}
\def \bex{\begin{exa}\label}
\def \eex{\end{exa}}
\def \bl{\begin{lem}\label}
\def \el{\end{lem}}
\def \bt{\begin{thm}\label}
\def \et{\end{thm}}
\def \bp{\begin{prop}\label}
\def \ep{\end{prop}}
\def \br{\begin{rem}\label}
\def \er{\end{rem}}
\def \bc{\begin{coro}\label}
\def \ec{\end{coro}}
\def \bd{\begin{de}\label}
\def \ed{\end{de}}

\newtheorem{thm}{Theorem}[section]
\newtheorem{prop}[thm]{Proposition}
\newtheorem{coro}[thm]{Corollary}
\newtheorem{conj}[thm]{Conjecture}
\newtheorem{exa}[thm]{Example}
\newtheorem{lem}[thm]{Lemma}
\newtheorem{rem}[thm]{Remark}
\newtheorem{de}[thm]{Definition}
\newtheorem{hy}[thm]{Hypothesis}
\makeatletter
\@addtoreset{equation}{section}
\def\theequation{\thesection.\arabic{equation}}
\makeatother
\makeatletter

\begin{center}{\Large \bf Axiomatic $G_{1}$-vertex algebras}
\end{center}

\begin{center}{Haisheng Li\footnote{Partially supported by NSF grant
DMS-9970496 and a grant from Rutgers Research Council}\\
Department of Mathematical Sciences, Rutgers University, Camden, NJ 08102\\
and\\
Department of Mathematics, Harbin Normal University, Harbin, China}
\end{center}

\begin{abstract}
Inspired by
Borcherds' work on ``$G$-vertex algebras,'' we formulate and 
study an axiomatic counterpart of Borcherds' notion of 
$G$-vertex algebra for the simplest nontrivial 
elementary vertex group, which we denote by $G_{1}$. 
Specifically, we formulate 
a notion of axiomatic $G_{1}$-vertex algebra, prove 
certain basic properties and give certain examples.
The notion of axiomatic $G_{1}$-vertex algebra is 
a nonlocal generalization of the notion of vertex algebra.
We also show how to construct axiomatic $G_{1}$-vertex
algebras from a set of compatible $G_{1}$-vertex operators.

The results of this paper had been reported in June 2001, at the
International Conference on Lie Algebras in the Morningside center,
Beijing, China, and been reported on November 30, 2001, in
the Quantum Mathematics Seminar, at Rutgers-New Brunswick.
We noticed that a paper of Bakalov and Kac appeared
today (math.QA/0204282) on noncommutative generalizations
of vertex algebras, which has certain overlaps with the current
paper. On the other hand, most of their results are othorgonal
to the results of this paper.
\end{abstract}

\section{Introduction}
It has been well known that vertex (operator) algebras introduced in
[B1] and [FLM]
are mathematical counterparts of chiral algebras in 2-dimensional 
quantum conformal field theory (cf. [BPZ]). 
Later, higher dimensional analogues of vertex algebras (or chiral algebras),
which are expected to play the same role in higher dimensional 
quantum field theory as vertex algebras play
in 2-dimensional quantum field theory,
were also established by Borcherds in [B2] 
by introducing a notion of $G$-vertex algebra. 
In this notion, $G$ is what Borcherds called an elementary vertex group
and a $G$-vertex algebra is 
an ``associative algebra'' in a certain ``relaxed multilinear category''
associated to $G$.
For the simplest nontrivial 
elementary vertex group $G$, which is denoted here by $G_{1}$, 
as it was proved in [Sn] (cf. [B2]), the notion of
{\em commutative} $G_{1}$-vertex algebra is equivalent to 
the notion of ordinary vertex algebra.
On the other hand, it had been known (earlier) that 
vertex (operator) algebras are analogous to
commutative associative algebras, as the Jacobi identity
for ordinary vertex algebras amounts to 
certain commutativity and associativity properties
(see [FLM], [FHL], [DL], [Li1]).
(Ordinary vertex algebras are also analogous to
Lie algebras in many aspects.)
In view of this analogy, a natural exercise is to establish
the corresponding analogues of
noncommutative associative algebras, or more or less,
to establish the axiomatic analogues of
Borcherds' $G_{1}$-vertex algebras.

In this paper we formulate and study
a notion of what we call axiomatic $G_{1}$-vertex algebra,
where axiomatic $G_{1}$-vertex algebras are
the corresponding analogues of noncommutative associative algebras
in contrast with the analogy between vertex algebras and 
commutative associative algebras.
This notion is defined by using all the axioms in the definition of the
notion of vertex algebra except for the Jacobi identity which is
replaced by the {\em weak associativity}:
For any algebra elements $u,v,w$, there exists a nonnegative integer
$l$ (depending only on $u$ and $w$) such that
\begin{eqnarray}
(x_{0}+x_{2})^{l}Y(u,x_{0}+x_{2})Y(v,x_{2})w=
(x_{0}+x_{2})^{l}Y(Y(u,x_{0})v,x_{2})w.
\end{eqnarray}
(See [DL] and [Li1] for this property.)
%a version of the physicists' operator product expansion formula.
It is expected that this notion of axiomatic $G_{1}$-vertex algebra 
is equivalent to the Borcherds' notion 
of $G_{1}$-vertex algebra. 

In terms of this notion, ordinary vertex algebras are exactly 
axiomatic $G_{1}$-vertex algebras
that also satisfy a certain weak commutativity property,
which was discovered to be an axiom in [DL] and [Li1].
Such a weak commutativity is also called {\em locality} in the literature.
%Ordinary vertex (super-)algebras (cf. [T]) 
%are also axiomatic $G_{1}$-vertex algebras.
Trivial examples of (nonlocal) axiomatic $G_{1}$-vertex algebras are
noncommutative associative algebras (with identity element).
We also give three constructions of (nonlocal) axiomatic
$G_{1}$-vertex algebras from ordinary vertex algebras.
All the three constructions are natural analogues of those in the
classical associative algebra theory.
In the classical theory, if $G$ is an abelian group,
any (commutative) $G$-graded associative algebra $A$ can be made a
(noncommutative) associative algebra by using a normalized $2$-cocycle
on $G$. Our first construction is an exact analogue
of this for a $G$-graded ordinary vertex algebra $V$
and a normalized $2$-cocycle $\epsilon$ on $G$.
Also, in the classical theory, if $A$ is a (commutative) associative
algebra and $n$ is a positive integer, we have the noncommutative 
matrix algebra $M(n, A)$ of all $n\times n$ matrices over $A$ for $n\ge 2$.
As our second construction, we show that for any ordinary vertex
algebra $V$, the vector space $M(n,V)$ has a natural axiomatic
$G_{1}$-vertex algebra structure. In fact, this follows from a 
general result.
Just as in [FHL] for ordinary vertex algebras, it can be easily shown 
that tensor products of axiomatic $G_{1}$-vertex algebras are
also axiomatic $G_{1}$-vertex algebras, see also [B2].
In particular, the tensor product of an ordinary vertex algebra
with an associative algebra is an axiomatic $G_{1}$-vertex algebra.
The axiomatic $G_{1}$-vertex algebra
$M(n,V)$ is naturally isomorphic to the tensor product 
axiomatic $G_{1}$-vertex algebra $V\otimes M(n,\C)$.
Now, let $A$ be an associative algebra acted by a group $G$ by automorphisms.
Associated to $A$ and $G$ there is an associative algebra 
called the cross product of $A$ with $G$ (or the skew algebra),
whose underlying vector space is $A\otimes \C[G]$.
Our third construction is an analogous cross product (or skew
product) construction of axiomatic $G_{1}$-vertex algebras from 
an ordinary vertex algebra acted by a group $G$.
For cross product axiomatic $G_{1}$-vertex algebras, 
we also derive a Jacobi-like identity and motivated by this, 
we define a notion of restricted (weak) axiomatic $G_{1}$-vertex algebra
by using a Jacobi-like identity as its main axiom.
This notion turns out to unify all the examples mentioned above.

There is a viewpoint about vertex (operator) algebras
which is that vertex (operator) algebras are
``algebras'' of vertex operators just as classical associative
algebras are algebras of linear operators. 
{}From this point of view, vertex operators
ought to give rise to vertex (operator) algebras 
just as linear operators naturally give rise to 
classical associative algebras.
In [Li1], for any abstract vector space $W$,
a notion of (weak) vertex operator on $W$ 
was defined and it was proved that 
any set of ``pairwise mutually local''
(weak) vertex operators on $W$ in a certain canonical way 
generates a vertex algebra with $W$ as a natural module.
This is an analogue of the classical fact that any set of 
pairwise commuting linear operators on $W$ generates a 
commutative associative algebra with $W$ as a module.
(See [Li2-3], [GL] for generalizations in certain directions.)
In the context of $G$-vertex algebras,
a theorem of Borcherds ([B2], Theorem 7.9) states that
any compatible set of vertex operators in a certain sense
on a vector space $W$ generates a $G$-vertex algebra acting on $W$.
Borcherds' theorem is also in the same spirit of the corresponding
theorems of [Li1], [Li2] and [GL], and it 
can be viewed as a noncommutative version of 
those corresponding theorems.
In this paper, we also study (weak) vertex operators 
as defined in [Li1], but they are renamed as (weak) $G_{1}$-vertex 
operators according to Borcherds' notion of $G$-vertex algebra.
We first define a notion of compatibility,
where the notion of compatibility is 
weaker than the notion of locality and it
asserts that the operator product expansion is of a certain form.
We then prove an analogous theorem (Theorem \ref{tgeneratingthem})
of Borcherds'. To prove this, we prove that
any closed compatible space, in a certain sense, of (weak) 
$G_{1}$-vertex operators on a vector space $W$ has a natural 
(weak) axiomatic $G_{1}$-vertex algebra 
structure with $W$ as a natural module (Theorem \ref{tclosed}),
which is analogous to a result obtained in [Li1] (cf. [MN])
for ordinary vertex algebras.

%In classical case, the composition or product of two linear 
%operators is a natural linear operator. For two vertex operators, 
%the existence of product expansion in a certain form
%is an assumption while locality is a sufficient condition.
%It is a trivial fact that {\em any} set of 
%linear operators on $W$ generates a 
%associative algebra with $W$ as a module.
%However, an {\em arbitrary} set of $G_{1}$-vertex operators
%does not ``generate'' an axiomatic $G_{1}$-vertex algebra.

In [EK], Etingof and Kazhdan established and studied a notion of 
quantum vertex operator algebra where it was proved that
a certain $h$-adic (topological) version of the weak associativity 
holds. To a certain extent, quantum vertex operator algebras 
are $h$-adic (topological) axiomatic $G_{1}$-vertex algebras.
Much of the current work can be carried on to quantum vertex operator
algebras and details will appear in a coming paper.

Recently, there has been active research in physics 
on noncommutative field theory 
(field theory on noncommutative manifolds) (cf. [DN]). 
It seems that noncommutative field theories
are related to (nonlocal) axiomatic $G_{1}$-vertex algebras.
This is also part of our motivation systematically to 
study axiomatic $G_{1}$-vertex algebras.

This paper is organized in the following manner.
In Section 2, we define a notion of axiomatic $G_{1}$-vertex algebra
and present certain basic properties. 
In Section 3, we discuss various examples
and we introduce a notions of restricted (weak) 
axiomatic $G_{1}$-vertex algebra
to unify many examples. In Section 4, we define the notion of module
and present certain basic properties. In Section 5, we show how 
to construct axiomatic $G_{1}$-vertex algebras from a set of
compatible $G_{1}$-vertex operators on a vector space $W$ 
and prove the main results.

\section{Axiomatic $G_{1}$-vertex algebras}
In this section, we define the notion of (weak) axiomatic 
$G_{1}$-vertex algebra and we establish certain basic properties 
analogous to those (cf. [DL], [FHL], [LL], [Li1]) 
for ordinary vertex algebras.

Let $x, y, z, x_{i}, y_{i}, z_{i}$, $i=0,1,\dots$ be
mutually commuting independent formal variables throughout this paper.
We shall use the standard formal variable notations and conventions
as defined in [FLM] and [FHL].
In particular, for a vector space $U$, 
\begin{eqnarray}
U[[x,x^{-1}]]&=&\left\{ \sum_{n\in \Z}a(n)x^{n}\;|\; a(n)\in U\right\},\\
U((x))&=&\left\{ \sum_{n\ge r}a(n)x^{n}\;|\; r\in \Z,\; a(n)\in
U\right\}
\subset U[[x,x^{-1}]],\\
U[[x]]&=&\left\{ \sum_{n\ge 0}a(n)x^{n}\;|\; a(n)\in U\right\}\subset U((x)).
\end{eqnarray}
The spaces $U[[x_{1},x_{1}^{-1},\dots,x_{n},x_{n}^{-1}]]$,
$U((x_{1},\dots,x_{n}))$ and $U[[x_{1},\dots,x_{n}]]$ are 
also defined in the obvious ways.

The formal delta function is the following formal series:
\begin{eqnarray}	
\delta(x)=\sum_{n\in \Z}x^{n}\in \C[[x,x^{-1}]].
\end{eqnarray}
In formal calculus the following {\em binomial expansion convention} 
is implemented:
\begin{eqnarray}
(x_{1}\pm x_{2})^{n}=\sum_{i\ge 0}{n\choose i}(\pm 1)^{i}x_{1}^{n-i}x_{2}^{i}
\in \C[[x_{1},x_{2},x_{1}^{-1}]],
\end{eqnarray}
where ${n\choose i}={1\over i!}n(n-1)\cdots (n+1-i)$ for $n\in \Z,\; i\in \N$.
Furthermore, by definition
\begin{eqnarray}
\delta\left(\frac{x_{1}-x_{2}}{x_{0}}\right)
=\sum_{n\in \Z}x_{0}^{-n}(x_{1}-x_{2})^{n}
=\sum_{n\in \Z}\sum_{i\ge 0}{n\choose i}(-1)^{i}x_{0}^{-n}x_{1}^{n-i}x_{2}^{i}.
\end{eqnarray}
Then we have 
\begin{eqnarray}
x_{0}^{-1}\delta\left(\frac{x_{1}-x_{2}}{x_{0}}\right)
-x_{0}^{-1}\delta\left(\frac{x_{2}-x_{1}}{-x_{0}}\right)
=x_{2}^{-1}\delta\left(\frac{x_{1}-x_{0}}{x_{2}}\right).
\end{eqnarray}

As it was emphasized in [FLM], in formal calculus,
associativity in general {\em does not hold} for products, and 
on the other hand, associativity {\em does hold} under 
the assumption that the products and their involved 
subproducts exist. For example, for three formal series $A, B$
and $C$, we have $A(BC)=(AB)C\;(=ABC)$ if $ABC$, $AB$ and $BC$ all exist.

The following is a reformulation of Proposition 3.4.2 of
[LL] with a slightly different proof
(cf. [Li1], part 3 of the proof of Proposition 2.2.4):

\bl{lformaljacobiidentity}
Let $U$ be a vector space and let
\begin{eqnarray}
A(x_{1},x_{2})&\in& U((x_{1}))((x_{2})),\\
B(x_{1},x_{2})&\in& U((x_{2}))((x_{1})),\\
C(x_{0},x_{2})&\in& U((x_{2}))((x_{0})).
\end{eqnarray}
Then 
\begin{eqnarray}\label{eformaljacobiABC}
x_{0}^{-1}\delta\left(\frac{x_{1}-x_{2}}{x_{0}}\right)A(x_{1},x_{2})
-x_{0}^{-1}\delta\left(\frac{x_{2}-x_{1}}{-x_{0}}\right)B(x_{1},x_{2})
=x_{2}^{-1}\delta\left(\frac{x_{1}-x_{0}}{x_{2}}\right)C(x_{0},x_{2})
\end{eqnarray}
if and only if there exist nonnegative integers 
$k$ and $l$ such that
\begin{eqnarray}
(x_{1}-x_{2})^{k}A(x_{1},x_{2})&=&(x_{1}-x_{2})^{k}B(x_{1},x_{2}),
\label{eformalA=B}\\
(x_{0}+x_{2})^{l}A(x_{0}+x_{2},x_{2})&=&(x_{0}+x_{2})^{l}C(x_{0},x_{2}).
\label{eformalA=C}
\end{eqnarray}
\el

{\bf Proof.} Let $k$ and $l$ be nonnegative integers such that
$$x_{0}^{k}C(x_{0},x_{2})\in U((x_{2}))[[x_{0}]],\;\;\;\;
x_{1}^{l}B(x_{1},x_{2})\in U((x_{2}))[[x_{1}]].$$
Then (\ref{eformalA=B}) follows from
(\ref{eformaljacobiABC}) by applying $\Res_{x_{0}}x_{0}^{k}$
and (\ref{eformalA=C}) follows from
(\ref{eformaljacobiABC}) by applying $\Res_{x_{1}}x_{1}^{l}$. 
Now we shall show that
(\ref{eformaljacobiABC}) also follows from (\ref{eformalA=B}) and 
(\ref{eformalA=C}). Let $l'$ be a nonnegative integer such that $l'\ge l$
and 
$x_{1}^{l'}B(x_{1},x_{2})\in U((x_{2}))[[x_{1}]]$.
Then (\ref{eformalA=B}) and (\ref{eformalA=C}) with $l$ being replaced
by $l'$ still hold, and furthermore,
$$(x_{1}-x_{2})^{k}x_{1}^{l'}A(x_{1},x_{2})
=(x_{1}-x_{2})^{k}x_{1}^{l'}B(x_{1},x_{2})\in U((x_{2}))[[x_{1}]]\cap
U((x_{1}))((x_{2})).$$
Thus
\begin{eqnarray}\label{echangesub}
\left[(x_{1}-x_{2})^{k}x_{1}^{l'}A(x_{1},x_{2})\right]|_{x_{1}=x_{2}+x_{0}}
=\left[(x_{1}-x_{2})^{k}x_{1}^{l'}A(x_{1},x_{2})\right]|_{x_{1}=x_{0}+x_{2}}.
\end{eqnarray}
Using (\ref{eformalA=B}), (\ref{echangesub}) and (\ref{eformalA=C}) we get
\begin{eqnarray}
& &x_{0}^{-1}\delta\left(\frac{x_{1}-x_{2}}{x_{0}}\right)
x_{0}^{k}x_{1}^{l'}A(x_{1},x_{2})
-x_{0}^{-1}\delta\left(\frac{x_{2}-x_{1}}{-x_{0}}\right)
x_{0}^{k}x_{1'}^{l}B(x_{1},x_{2})\nonumber\\
&=&x_{0}^{-1}\delta\left(\frac{x_{1}-x_{2}}{x_{0}}\right)
(x_{1}-x_{2})^{k}x_{1}^{l'}A(x_{1},x_{2})
-x_{0}^{-1}\delta\left(\frac{x_{2}-x_{1}}{-x_{0}}\right)
(x_{1}-x_{2})^{k}x_{1}^{l'}B(x_{1},x_{2})\nonumber\\
&=&x_{2}^{-1}\delta\left(\frac{x_{1}-x_{0}}{x_{2}}\right)
\left[ (x_{1}-x_{2})^{k}x_{1}^{l'}A(x_{1},x_{2})\right]\nonumber\\
&=&x_{2}^{-1}\delta\left(\frac{x_{1}-x_{0}}{x_{2}}\right)
\left[(x_{1}-x_{2})^{k}x_{1}^{l'}A(x_{1},x_{2})\right]|_{x_{1}=x_{2}+x_{0}}
\nonumber\\
&=&x_{2}^{-1}\delta\left(\frac{x_{1}-x_{0}}{x_{2}}\right)
\left[(x_{1}-x_{2})^{k}x_{1}^{l'}A(x_{1},x_{2})\right]|_{x_{1}=x_{0}+x_{2}}
\nonumber\\
&=&x_{2}^{-1}\delta\left(\frac{x_{1}-x_{0}}{x_{2}}\right)
\left[x_{0}^{k}(x_{0}+x_{2})^{l'}A(x_{0}+x_{2},x_{2})\right]
\nonumber\\
&=&x_{2}^{-1}\delta\left(\frac{x_{1}-x_{0}}{x_{2}}\right)
\left[x_{0}^{k}(x_{0}+x_{2})^{l'}C(x_{0},x_{2})\right]
\nonumber\\
&=&x_{2}^{-1}\delta\left(\frac{x_{1}-x_{0}}{x_{2}}\right)
x_{0}^{k}x_{1}^{l'}C(x_{0},x_{2}).
\end{eqnarray}
Multiplying both sides by $x_{0}^{-k}x_{1}^{-l'}$ we obtain
(\ref{eformaljacobiABC}).
$\;\;\;\;\Box$

The notion of axiomatic $G_{1}$-vertex algebra is defined as follows:

\bd{dnoncomva}
{\em An {\em axiomatic $G_{1}$-vertex algebra} is 
a vector space $V$ equipped with a linear map 
\begin{eqnarray}
Y(\cdot,x):V&\rightarrow& (\End V)[[x,x^{-1}]]\nonumber\\
v&\mapsto& Y(v,x)=\sum_{n\in \Z}v_{n}x^{-n-1}
\end{eqnarray}
and equipped with a distinguished vector ${\bf 1}\in V$ 
such that all the following axioms hold: For $u,v\in V$, 
\begin{eqnarray}
u_{n}v=0\;\;\;\mbox{  for $n$ sufficiently large};
\end{eqnarray}
\begin{eqnarray}
Y({\bf 1},x)=1;
\end{eqnarray}
for $v\in V$,
\begin{eqnarray}\label{edefcreation}
Y(v,x){\bf 1}\in V[[x]]\;\;\;\mbox{ and }\;\;\;
\lim_{x\rightarrow 0}Y(v,x){\bf 1}(=v_{-1}{\bf 1})=v;
\end{eqnarray}
and for $u,w\in V$, there exists $l\in \N$ such that for all $v\in V$,
\begin{eqnarray}\label{eweakassoc}
(x_{0}+x_{2})^{l}Y(u,x_{0}+x_{2})Y(v,x_{2})w
=(x_{0}+x_{2})^{l}Y(Y(u,x_{0})v,x_{2})w
\end{eqnarray}
(the {\em weak associativity}).}
\ed

Note that the integer $l$ in (\ref{eweakassoc}) only depends on $u$
and $w$, but not $v$.
A {\em weak axiomatic $G_{1}$-vertex algebra} 
satisfies all the axioms for an axiomatic $G_{1}$-vertex algebra
except that the weak associativity axiom is replaced by
the weaker one: {\em For any $u,v, w\in V$ there exists 
a nonnegative integer $l$ such that (\ref{eweakassoc}) holds}.

Of course, a finite-dimensional weak axiomatic $G_{1}$-vertex algebra
is automatically an axiomatic $G_{1}$-vertex algebra.

\br{requivalenceandnotion}
{\em In the notion of (weak) axiomatic $G_{1}$-vertex algebra,
$G_{1}$ represents the simplest nontrivial
elementary vertex group defined in [B2]; roughly speaking,
it is the pair $(H_{1},K_{1})$, where $H_{1}=\C[D]$, being considered as
the universal enveloping algebra of the $1$-dimensional Lie algebra $\C D$,
is a cocommutative Hopf algebra and $K_{1}=\C[[x]][x^{-1}]=\C((x))$,
a commutative associative algebra and an $H_{1}$-module 
with $D$ acting as $d/dx$. It is expected that
the notion of axiomatic $G_{1}$-vertex algebra is equivalent
to Borcherds' notion of $G_{1}$-vertex algebra defined in [B2].
{\em Within this paper},
when there is no confusion, we shall take the liberty simply to use
the term ``$G_{1}$-vertex algebra.''}
\er

\br{rweakcommassoc}
{\em Recall from [B1] and [FLM] (cf. [Li1]) that
a {\em vertex algebra} is a vector space $V$ 
such that all the axioms for an axiomatic $G_{1}$-vertex algebra except 
for the weak associativity hold and such that for $u,v\in V$,
\begin{eqnarray}\label{ejacobidef}
& &x_{0}^{-1}\delta\left(\frac{x_{1}-x_{2}}{x_{0}}\right)
Y(u,x_{1})Y(v,x_{2})-x_{0}^{-1}\delta\left(\frac{x_{2}-x_{1}}{-x_{0}}\right)
Y(v,x_{2})Y(u,x_{1})\nonumber\\
&=&x_{2}^{-1}\delta\left(\frac{x_{1}-x_{0}}{x_{2}}\right)
Y(Y(u,x_{0})v,x_{2})
\end{eqnarray}
(the {\em Jacobi identity}).
In view of Lemma \ref{lformaljacobiidentity} (cf. [FHL], [DL], [Li1]),
the Jacobi identity (\ref{ejacobidef}) is equivalent to
the following {\em weak commutativity and weak associativity}:
For $u,v\in V$, there exists $k\in \N$ such that
\begin{eqnarray}
(x_{1}-x_{2})^{k}Y(u,x_{1})Y(v,x_{2})=(x_{1}-x_{2})^{k}Y(v,x_{2})Y(u,x_{1});
\end{eqnarray}
for $u, w\in V$, there exists $l\in \N$ such that for all $v\in V$,
\begin{eqnarray}
(x_{0}+x_{2})^{l}Y(u,x_{0}+x_{2})Y(v,x_{2})w
=(x_{0}+x_{2})^{l}Y(Y(u,x_{0})v,x_{2})w.
\end{eqnarray}
In view of this, ordinary vertex algebras 
are axiomatic $G_{1}$-vertex algebras. Furthermore,
ordinary vertex algebras are analogous to commutative 
associative algebras while axiomatic $G_{1}$-vertex algebras are 
analogous to associative algebras.}
\er

\bex{exampleassociativealgebra}
{\em Just as ordinary vertex algebras can be constructed from 
commutative associative algebras with identity element 
equipped with a derivation (cf. [B1]),
axiomatic $G_{1}$-vertex algebras can be constructed from
associative algebras with identity element equipped with a derivation.
Let $A$ be an associative algebra with identity element $1$
equipped with a derivation $d$ (possibly zero). Define a linear map
\begin{eqnarray}
Y(\cdot,x): A\rightarrow (\End A)[[x]]\subset (\End A)[[x,x^{-1}]]
\end{eqnarray}
by
\begin{eqnarray}\label{edefvertexoperatorexample}
Y(a,x)b=\left(e^{xd}a\right)b=\sum_{i\ge 0}{1\over i!}(d^{i}a)bx^{i}.
\end{eqnarray}
All the axioms except for the weak associativity clearly hold.
Since $d$ is a derivation of $A$, $e^{xd}$ is an automorphism 
of the associative algebra $A[[x]]$ (by considering $d$ as a
derivation of $A[[x]]$ with $d(x)=0$), so that for $a,b,c\in A$,
\begin{eqnarray}
Y(a,x_{0}+x_{2})Y(b,x_{2})c
&=&\left(e^{(x_{0}+x_{2})d}a\right)\left(e^{x_{2}d}b\right)c\nonumber\\
&=&\left(e^{x_{2}d}\left(\left(e^{x_{0}d}a\right)b\right)\right)c\nonumber\\
&=&Y(Y(a,x_{0})b,x_{2})c.
\end{eqnarray}
This proves the weak associativity, so $A$ equipped with the distinguished vector $1$ and 
the linear map $Y$ defined in (\ref{edefvertexoperatorexample})
is an axiomatic $G_{1}$-vertex algebra.
In particular, by taking $d=0$ we see that any associative algebra
with identity is an axiomatic $G_{1}$-vertex algebra.}
\eex

Next, we give some consequences of the definition.
First, as  in [Li1] and [LL] for vertex algebras
we have the following $D$-bracket-derivative formula
(which in fact follows from the same proof of [LL]):

\bp{pcon1}
Let $V$ be a weak axiomatic $G_{1}$-vertex algebra. Define
a linear operator $D$ on $V$ by
\begin{eqnarray}
D(v)=v_{-2}{\bf 1}\left(=\left({d\over dx}Y(v,x){\bf 1}\right)|_{x=0}\right)
\;\;\;\mbox{ for }v\in V.
\end{eqnarray}
Then
\begin{eqnarray}\label{edproperty}
[D,Y(v,x)]=Y(D(v),x)={d\over dx}Y(v,x)\;\;\;\mbox{ for }v\in V.
\end{eqnarray}
\ep

{\bf Proof.}  Let $u,v\in V$. Then there exists $l\in \N$ such that
\begin{eqnarray}
(x_{2}+x_{0})^{l}Y(Y(u,x_{0}){\bf 1},x_{2})v
=(x_{2}+x_{0})^{l}Y(u,x_{0}+x_{2})Y({\bf 1},x_{2})v.
\end{eqnarray}
With $Y({\bf 1},x)=1$  we have
\begin{eqnarray}\label{e2.25}
(x_{2}+x_{0})^{l}Y(Y(u,x_{0}){\bf 1},x_{2})v
=(x_{0}+x_{2})^{l}Y(u,x_{0}+x_{2})v.
\end{eqnarray}
We may assume that $x^{l}Y(u,x)v\in V[[x]]$ by replacing $l$
with a bigger integer if necessary, so that
\begin{eqnarray}
(x_{0}+x_{2})^{l}Y(u,x_{0}+x_{2})v=(x_{2}+x_{0})^{l}Y(u,x_{2}+x_{0})v.
\end{eqnarray}
Then (\ref{e2.25}) can be also written as
\begin{eqnarray}
(x_{2}+x_{0})^{l}Y(Y(u,x_{0}){\bf 1},x_{2})v
=(x_{2}+x_{0})^{l}Y(u,x_{2}+x_{0})v.
\end{eqnarray}
Multiplying both sides by $(x_{2}+x_{0})^{-l}$ we get
\begin{eqnarray}
Y(Y(u,x_{0}){\bf 1},x_{2})v
=Y(u,x_{2}+x_{0})v=e^{x_{0}{d\over dx_{2}}}Y(u,x_{2})v.
\end{eqnarray}
Extracting the coefficient of $x_{0}$ we obtain
\begin{eqnarray}
Y(D(u),x_{2})v={d\over dx_{2}}Y(u,x_{2})v.
\end{eqnarray}
This proves the second equality of (\ref{edproperty}).

For the first equality, let $u,v\in V$ and let $l\in \N$ be such that
\begin{eqnarray}
(x_{2}+x_{0})^{l}Y(Y(u,x_{0})v,x_{2}){\bf 1}=
(x_{2}+x_{0})^{l}Y(u,x_{0}+x_{2})Y(v,x_{2}){\bf 1}.
\end{eqnarray}
In view of the creation property, $Y(Y(u,x_{0})v,x_{2}){\bf 1}$ 
involves only nonnegative powers of $x_{2}$, so that we may multiply 
both sides by $(x_{0}+x_{2})^{-l}$ to get
\begin{eqnarray}
Y(Y(u,x_{0})v,x_{2}){\bf 1}=
Y(u,x_{0}+x_{2})Y(v,x_{2}){\bf 1}.
\end{eqnarray}
In view of the Taylor theorem we have
\begin{eqnarray}
Y(Y(u,x_{0})v,x_{2}){\bf 1}
=e^{x_{2}{\partial\over\partial x_{0}}}Y(u,x_{0})Y(v,x_{2}){\bf 1}.
\end{eqnarray}
Extracting the coefficient of $x_{2}$ from both sides 
and using the creation property we get
\begin{eqnarray}
D(Y(u,x_{0})v)
=Y(u,x_{0})D(v)+{d\over dx_{0}}Y(u,x_{0})v.
\end{eqnarray}
That is,
\begin{eqnarray}
[D,Y(u,x_{0})]v={d\over dx_{0}}Y(u,x_{0})v.
\end{eqnarray}
This proves the first equality of (\ref{edproperty}), completing the proof.
$\;\;\;\;\Box$

Combining Proposition \ref{pcon1} with the Taylor theorem 
we immediately have the first part of the following Proposition:

\bc{cdoperator}
Let $V$ be a weak axiomatic $G_{1}$-vertex algebra and 
let $D\in \End V$ be defined as in 
Proposition \ref{pcon1}. Then for $v\in V$,
\begin{eqnarray}
& &e^{xD}Y(v,x_{1})e^{-xD}=Y(e^{xD}v,x_{1})=Y(v,x_{1}+x),
\label{econjugationformula1}\\
& &Y(v,x){\bf 1}=e^{xD}v.\label{ecreationwithd}
\end{eqnarray}
\ec

{\pf} Applying the second equality of (\ref{econjugationformula1}) 
to ${\bf 1}$, and then setting $x_{1}=0$ and using the creation property
 we obtain (\ref{ecreationwithd}).
$\;\;\;\;\Box$

\br{rweakvertexalgebra}
{\em Recall from [LL] that a {\em weak vertex algebra} is a vector
space $V$ equipped with a linear map $Y$ from $V$ to $(\End V)[[x,x^{-1}]]$
and equipped with a distinguished vector ${\bf 1}$ such that
$Y({\bf 1},x)=1$ and such that (\ref{edefcreation}) and (\ref{edproperty})
hold. In view of Proposition \ref{pcon1},
any weak axiomatic $G_{1}$-vertex algebra 
is a weak vertex algebra.}
\er

Let $V$ be a weak axiomatic $G_{1}$-vertex algebra. A {\em subalgebra}
of $V$ is a subspace $U$ such that 
\begin{eqnarray}
& &{\bf 1}\in U,\\
& &u_{n}u'\in U\;\;\;\mbox{ for }u,u'\in U,\; n\in \Z.
\end{eqnarray}
Then $U$ itself equipped with the linear map $Y$ restricted to $U$ 
is a weak axiomatic $G_{1}$-vertex algebra.

Let $U$ be a subspace of $V$. We define the stabilizer 
${\rm Stab}(U)$ of $U$ in $V$ as
\begin{eqnarray}
{\rm Stab}(U)=\{ v\in V\;|\; v_{n}U\subset U\;\;\;\mbox{ for all }n\in \Z\}.
\end{eqnarray}
Then $U$ is a subalgebra if and only if ${\bf 1}\in U$ and 
$U\subset {\rm Stab}(U)$.

\bl{lstablizer}
The stabilizer ${\rm Stab}(U)$ of $U$ in $V$ is a subalgebra.
\el

{\bf Proof.}  Clearly, ${\bf 1}\in {\rm Stab}(U)$. 
Let $a,b\in {\rm Stab}(U)$ and $u\in U$.
Then there exists a nonnegative integer $l$ such that
\begin{eqnarray}
(x_{2}+x_{0})^{l}Y(Y(a,x_{0})b,x_{2})u=
(x_{2}+x_{0})^{l}Y(a,x_{0}+x_{2})Y(b,x_{2})u.
\end{eqnarray}
With $a,b\in {\rm Stab}(U)$, we have
\begin{eqnarray*}
(x_{2}+x_{0})^{l}Y(a,x_{0}+x_{2})Y(b,x_{2})u\in 
U[[x_{0},x_{0}^{-1},x_{2},x_{2}^{-1}]],
\end{eqnarray*}
so 
\begin{eqnarray}
(x_{2}+x_{0})^{l}Y(Y(a,x_{0})b,x_{2})u\in 
U[[x_{0},x_{0}^{-1},x_{2},x_{2}^{-1}]].
\end{eqnarray}
Then
\begin{eqnarray}
(x_{2}+x_{0})^{l}Y(Y(a,x_{0})b,x_{2})u\in 
V((x_{2}))((x_{0}))\cap U[[x_{0},x_{0}^{-1},x_{2},x_{2}^{-1}]]
=U((x_{2}))((x_{0})).
\end{eqnarray}
Consequently, 
$Y(Y(a,x_{0})b,x_{2})u\in U((x_{2}))((x_{0}))$, since
$$(x_{2}+x_{0})^{-l}F(x_{0},x_{2})\in U((x_{2}))((x_{0}))
\;\;\;\mbox{ for any }F(x_{0},x_{2})\in U((x_{2}))((x_{0})).$$
Then $(a_{m}b)_{n}u\in U$ for all $m,n\in \Z$.
Thus $a_{m}b\in {\rm Stab}(U)$ for all $m\in \Z$. 
Therefore, ${\rm Stab}(U)$ is a subalgebra.
$\;\;\;\;\;\Box$

Let $S$ be a subset of a weak axiomatic $G_{1}$-vertex algebra $V$.
Denote by $\<S\>$ the {\em subalgebra of $V$ generated by $S$}, 
which is by definition the smallest subalgebra of
$V$ containing $S$.

\bp{psubalgebragenerated}
For any subset $S$ of $V$, the subalgebra $\<S\>$ generated by $S$ 
is linearly spanned by vectors
\begin{eqnarray}\label{espannform}
u^{(1)}_{n_{1}}\cdots u^{(r)}_{n_{r}}{\bf 1}
\end{eqnarray}
for $r\ge 0,\; u^{(i)}\in S,\; n_{1},\dots,n_{r}\in \Z$.
\ep

{\bf Proof.} Let $U$ be the subspace linearly spanned by vectors in
(\ref{espannform}). Since, any subalgebra that contains $S$ must
contain $U$, we have $U\subset \<S\>$. To prove $\<S\>\subset U$,
since $S\subset U$, it suffices to prove that $U$ is a subalgebra.
Since $S\subset {\rm Stab}(U)$ and 
${\rm Stab}(U)$ is a subalgebra 
(Lemma \ref{lstablizer}), we have $\<S\>\subset {\rm Stab}(U)$.
Consequently, $U\subset \<S\>\subset {\rm Stab}(U)$. 
Then $U$ is a subalgebra (clearly, ${\bf 1}\in U$), so that
$\<S\>\subset U$. This proves $U=\<S\>$, 
completing the proof.
$\;\;\;\;\Box$

For ordinary vertex algebras, due to the Borcherds' commutator formula,
we know (cf. [FHL]) that vertex operators
$Y(u,x_{1})$ and $Y(v,x_{2})$ commute if and only if 
$u_{i}v=0$ for $i\ge 0$. For (weak) axiomatic $G_{1}$-vertex algebras,
we in general do not have Borcherds' commutator formula.
Nevertheless, here we have:

\bp{pcommutativity}
Let $V$ be a weak axiomatic $G_{1}$-vertex algebra and 
let $u,v\in V,\; k\in \N,\; q\in \C^{*}$. Then 
\begin{eqnarray}\label{equsilocalityprop}
(x_{1}-x_{2})^{k}Y(u,x_{1})Y(v,x_{2})=q
(x_{1}-x_{2})^{k}Y(v,x_{2})Y(u,x_{1})
\end{eqnarray}
if and only if 
\begin{eqnarray}
& &x^{k}Y(u,x)v\in V[[x]],\label{etruncationprop1}\\
& &Y(u,x)v=qe^{xD}Y(v,-x)u.\label{eskewsymmetry}
\end{eqnarray}
In particular, $[Y(u,x_{1}),Y(v,x_{2})]=0$ if and only if 
$Y(u,x)v\in V[[x]]$ and $Y(u,x)v=e^{xD}Y(v,-x)u$.
\ep

{\bf Proof.} Assume (\ref{equsilocalityprop}) holds. 
For any $w\in V$, in view of Lemma \ref{lformaljacobiidentity}, 
(\ref{equsilocalityprop}),
together with the weak associativity relation
\begin{eqnarray*}
(x_{0}+x_{2})^{l}Y(u,x_{0}+x_{2})Y(v,x_{2})w=
(x_{0}+x_{2})^{l}Y(Y(u,x_{0})v,x_{2})w
\end{eqnarray*}
for some nonnegative integer $l$, gives
\begin{eqnarray}\label{ejacobiuvw}
& &x_{0}^{-1}\delta\left(\frac{x_{1}-x_{2}}{x_{0}}\right)
Y(u,x_{1})Y(v,x_{2})w
-x_{0}^{-1}\delta\left(\frac{x_{2}-x_{1}}{-x_{0}}\right)
qY(v,x_{2})Y(u,x_{1})w\nonumber\\
&=&x_{2}^{-1}\delta\left(\frac{x_{1}-x_{0}}{x_{2}}\right)
Y(Y(u,x_{0})v,x_{2})w.
\end{eqnarray}
With (\ref{equsilocalityprop}), after multiplied 
by $(x_{1}-x_{2})^{k}$ the left-hand side of (\ref{ejacobiuvw})
involves only nonnegative powers of $x_{0}$, so is the right-hand side.
Then (by taking $\Res_{x_{1}}$)
\begin{eqnarray}
x_{0}^{k}Y(Y(u,x_{0})v,x_{2})w\in V((x_{2}))[[x_{0}]]
\;\;\;\mbox{ for }w\in V.
\end{eqnarray}
Since the vertex operator map $Y$ is injective
(from the creation property (\ref{edefcreation})),
we obtain (\ref{etruncationprop1}).
Similarly, (\ref{equsilocalityprop}),
together with the weak associativity relation
\begin{eqnarray*}
(-x_{0}+x_{1})^{l}qY(v,-x_{0}+x_{1})Y(u,x_{1})w=
(-x_{0}+x_{1})^{l}qY(Y(v,-x_{0})u,x_{1})w
\end{eqnarray*}
for some nonnegative integer $l$, gives
\begin{eqnarray}
& &x_{0}^{-1}\delta\left(\frac{x_{1}-x_{2}}{x_{0}}\right)
Y(u,x_{1})Y(v,x_{2})w
-x_{0}^{-1}\delta\left(\frac{x_{2}-x_{1}}{-x_{0}}\right)
qY(v,x_{2})Y(u,x_{1})w\nonumber\\
&=&x_{1}^{-1}\delta\left(\frac{x_{2}+x_{0}}{x_{1}}\right)
qY(Y(v,-x_{0})u,x_{1})w.
\end{eqnarray}
Using (\ref{econjugationformula1}) we get
\begin{eqnarray}
& &x_{1}^{-1}\delta\left(\frac{x_{2}+x_{0}}{x_{1}}\right)
qY(Y(v,-x_{0})u,x_{1})w\nonumber\\
&=&x_{1}^{-1}\delta\left(\frac{x_{2}+x_{0}}{x_{1}}\right)
qY(Y(v,-x_{0})u,x_{2}+x_{0})w\nonumber\\
&=&x_{1}^{-1}\delta\left(\frac{x_{2}+x_{0}}{x_{1}}\right)
qY(e^{x_{0}D}Y(v,-x_{0})u,x_{2})w.
\end{eqnarray}
Thus
\begin{eqnarray}
x_{2}^{-1}\delta\left(\frac{x_{1}-x_{0}}{x_{2}}\right)
Y(Y(u,x_{0})v,x_{2})w
=x_{1}^{-1}\delta\left(\frac{x_{2}+x_{0}}{x_{1}}\right)
qY(e^{x_{0}D}Y(v,-x_{0})u,x_{2})w,
\end{eqnarray}
which (by taking $\Res_{x_{1}}$) gives
\begin{eqnarray}
Y(Y(u,x_{0})v,x_{2})w=qY(e^{x_{0}D}Y(v,-x_{0})u,x_{2})w.
\end{eqnarray}
Again, with $Y$ being injective,
the skew-symmetry relation (\ref{eskewsymmetry}) follows immediately.

On the other hand, assume that (\ref{etruncationprop1}) and 
(\ref{eskewsymmetry}) hold. 
Let $w\in V$. There exists 
a nonnegative integer $l$ such that
$x^{l}Y(v,x)w\in V[[x]]$ and 
\begin{eqnarray}
(x_{0}+x_{2})^{l}Y(u,x_{0}+x_{2})Y(v,x_{2})w
&=&(x_{0}+x_{2})^{l}Y(Y(u,x_{0})v,x_{2})w,\label{eweakassocrelation1}\\
(-x_{0}+x_{1})^{l}Y(v,-x_{0}+x_{1})Y(u,x_{1})w
&=&(-x_{0}+x_{1})^{l}Y(Y(v,-x_{0})u,x_{1})w.\label{eweakassocrelation2}
\end{eqnarray}
Then 
\begin{eqnarray*}
x_{2}^{l}(x_{0}+x_{2})^{l}Y(Y(u,x_{0})v,x_{2})w
=x_{2}^{l}(x_{0}+x_{2})^{l}Y(u,x_{0}+x_{2})Y(v,x_{2})w
\in V((x_{0}))[[x_{2}]],
\end{eqnarray*}
hence
\begin{eqnarray}\label{eevaluationiden}
\left[x_{2}^{l}(x_{0}+x_{2})^{l}Y(Y(u,x_{0})v,x_{2})w\right]|_{x_{2}
=-x_{0}+x_{1}}
=\left[x_{2}^{l}(x_{0}+x_{2})^{l}Y(Y(u,x_{0})v,x_{2})w\right]|_{x_{2}
=x_{1}-x_{1}}.
\end{eqnarray}
{}From (\ref{eskewsymmetry}) and (\ref{econjugationformula1})
we have
\begin{eqnarray}\label{eqskewuvw}
qY(Y(v,-x_{0})u,x_{1})=Y(e^{-x_{0}D}Y(u,x_{0})v,x_{1})
=Y(Y(u,x_{0})v,x_{1}-x_{0}).
\end{eqnarray}
Using  (\ref{eweakassocrelation1}), (\ref{eweakassocrelation2}),
(\ref{eqskewuvw}) and (\ref{eevaluationiden}) we get
\begin{eqnarray}
& &x_{0}^{-1}\delta\left(\frac{x_{1}-x_{2}}{x_{0}}\right)
x_{0}^{k}x_{1}^{l}x_{2}^{l}Y(u,x_{1})Y(v,x_{2})w\nonumber\\
& &-x_{0}^{-1}\delta\left(\frac{x_{2}-x_{1}}{-x_{0}}\right)
x_{0}^{k}x_{1}^{l}x_{2}^{l}qY(v,x_{2})Y(u,x_{1})w\nonumber\\
&=&x_{0}^{-1}\delta\left(\frac{x_{1}-x_{2}}{x_{0}}\right)
x_{0}^{k}(x_{0}+x_{2})^{l}x_{2}^{l}Y(u,x_{0}+x_{2})Y(v,x_{2})w\nonumber\\
& &-x_{0}^{-1}\delta\left(\frac{x_{2}-x_{1}}{-x_{0}}\right)
x_{0}^{k}x_{1}^{l}(-x_{0}+x_{1})^{l}qY(v,-x_{0}+x_{1})Y(u,x_{1})w\nonumber\\
&=&x_{0}^{-1}\delta\left(\frac{x_{1}-x_{2}}{x_{0}}\right)
\left[x_{0}^{k}(x_{0}+x_{2})^{l}x_{2}^{l}Y(Y(u,x_{0})v,x_{2})w\right]
\nonumber\\
& &-x_{0}^{-1}\delta\left(\frac{x_{2}-x_{1}}{-x_{0}}\right)
\left[x_{0}^{k}x_{1}^{l}(-x_{0}+x_{1})^{l}qY(Y(v,-x_{0})u,x_{1})w\right]
\nonumber\\
&=&x_{0}^{-1}\delta\left(\frac{x_{1}-x_{2}}{x_{0}}\right)
\left[x_{0}^{k}(x_{0}+x_{2})^{l}x_{2}^{l}Y(Y(u,x_{0})v,x_{2})w\right]
\nonumber\\
& &-x_{0}^{-1}\delta\left(\frac{x_{2}-x_{1}}{-x_{0}}\right)
\left[x_{0}^{k}x_{1}^{l}(-x_{0}+x_{1})^{l}Y(Y(u,x_{0})v,x_{1}-x_{0})w
\right]\nonumber\\
&=&x_{2}^{-1}\delta\left(\frac{x_{1}-x_{0}}{x_{2}}\right)
\left[x_{0}^{k}(x_{0}+x_{2})^{l}x_{2}^{l}Y(Y(u,x_{0})v,x_{2})w\right]
\nonumber\\
& &+x_{0}^{-1}\delta\left(\frac{x_{2}-x_{1}}{-x_{0}}\right)
\left[x_{0}^{k}(x_{0}+x_{2})^{l}x_{2}^{l}Y(Y(u,x_{0})v,x_{2})w\right]
\nonumber\\
& &-x_{0}^{-1}\delta\left(\frac{x_{2}-x_{1}}{-x_{0}}\right)
\left[x_{0}^{k}x_{1}^{l}(-x_{0}+x_{1})^{l}Y(Y(u,x_{0})v,x_{1}-x_{0})w
\right]\nonumber\\
&=&x_{2}^{-1}\delta\left(\frac{x_{1}-x_{0}}{x_{2}}\right)
x_{0}^{k}(x_{0}+x_{2})^{l}x_{2}^{l}Y(Y(u,x_{0})v,x_{2})w\nonumber\\
& &+x_{0}^{-1}\delta\left(\frac{x_{2}-x_{1}}{-x_{0}}\right)
\left[x_{0}^{k}(x_{0}+x_{2})^{l}x_{2}^{l}Y(Y(u,x_{0})v,x_{2})w\right]
|_{x_{2}=-x_{0}+x_{1}}\nonumber\\
& &-x_{0}^{-1}\delta\left(\frac{x_{2}-x_{1}}{-x_{0}}\right)
\left[x_{0}^{k}x_{1}^{l}(-x_{0}+x_{1})^{l}Y(Y(u,x_{0})v,x_{1}-x_{0})w
\right]\nonumber\\
&=&x_{2}^{-1}\delta\left(\frac{x_{1}-x_{0}}{x_{2}}\right)
x_{0}^{k}x_{1}^{l}x_{2}^{l}Y(Y(u,x_{0})v,x_{2})w
\nonumber\\
& &+x_{0}^{-1}\delta\left(\frac{x_{2}-x_{1}}{-x_{0}}\right)
\left[x_{0}^{k}(x_{0}+x_{2})^{l}x_{2}^{l}Y(Y(u,x_{0})v,x_{2})w\right]
|_{x_{2}=x_{1}-x_{0}}\nonumber\\
& &-x_{0}^{-1}\delta\left(\frac{x_{2}-x_{1}}{-x_{0}}\right)
\left[x_{0}^{k}x_{1}^{l}(-x_{0}+x_{1})^{l}Y(Y(u,x_{0})v,x_{1}-x_{0})w
\right]\nonumber\\
&=&x_{2}^{-1}\delta\left(\frac{x_{1}-x_{0}}{x_{2}}\right)
x_{0}^{k}x_{1}^{l}x_{2}^{l}Y(Y(u,x_{0})v,x_{2})w.
\end{eqnarray}
Taking $\Res_{x_{0}}$, then using (\ref{etruncationprop1})
we obtain (\ref{equsilocalityprop}).
$\;\;\;\;\Box$

\br{rjacobiidentity2}
{\em As we have seen in the proof of Proposition \ref{pcommutativity},
for weak axiomatic $G_{1}$-vertex algebras,
the weak commutativity relation
(\ref{equsilocalityprop}) amounts to the following Jacobi identity
\begin{eqnarray}\label{ejacobiyuyvw}
& &x_{0}^{-1}\delta\left(\frac{x_{1}-x_{2}}{x_{0}}\right)Y(u,x_{1})Y(v,x_{2})w
-x_{0}^{-1}\delta\left(\frac{x_{2}-x_{1}}{-x_{0}}\right)q Y(v,x_{2})Y(u,x_{1})w
\nonumber\\
&=&x_{2}^{-1}\delta\left(\frac{x_{1}-x_{0}}{x_{2}}\right)Y(Y(u,x_{0})v,x_{2})w.
\end{eqnarray}}
\er

Notice that for any $u,v\in V$, there always exists  
a nonnegative integer $k$ such that $x^{k}Y(u,x)v\in V[[x]]$. Then
as an immediate consequence of Proposition \ref{pcommutativity} we have:

\bc{clocalityskewsymmetry}
Let $u,v\in V$. Then there exists a nonnegative integer $k$ such that
\begin{eqnarray}
(x_{1}-x_{2})^{k}[Y(u,x_{1}),Y(v,x_{2})]=0
\end{eqnarray}
if and only if 
\begin{eqnarray}
Y(u,x)v=e^{xD}Y(v,-x)u.\;\;\;\;\Box
\end{eqnarray}
\ec

\br{rLLsame}
{\em It was known (cf. [FHL], [Li1], [LL]) that in the theory of
ordinary vertex algebras, under the skew-symmetry,
the weak commutativity is equivalent to the weak associativity.
Proposition \ref{pcommutativity} and Corollary
\ref{clocalityskewsymmetry} are in the same spirit.}
\er

Let $S$ be a subset of a weak axiomatic $G_{1}$-vertex algebra $V$. 
We define the {\em localizer $L_{V}(S)$ of $S$ in $V$}
to consist of $v\in V$ such that for every $w\in S$ there exists a
nonnegative integer $k$ such that
$$(x_{1}-x_{2})^{k}[Y(v,x_{1}),Y(w,x_{2})]=0.$$
In view of Corollary \ref{clocalityskewsymmetry} we have
\begin{eqnarray}
L_{V}(S)=\{ v\in V\;|\; Y(v,x)w=e^{xD}Y(w,-x)v\;\;\mbox{ for every }w\in S\}.
\end{eqnarray}

\bp{plocalizer}
For any subset $S$ of $V$, the localizer $L_{V}(S)$ is a subalgebra.
\ep

{\bf Proof.} Clearly $1\in L_{V}(S)$, so we must prove that 
$u_{n}v\in L_{V}(S)$ for $u,v\in L_{V}(S),\; n\in \Z$.
In view of Proposition \ref{pcommutativity}, we must show
\begin{eqnarray}
Y(u_{n}v,x)w=e^{xD}Y(w,-x)u_{n}v\;\;\;\mbox{ for }w\in S.
\end{eqnarray}

Let $u,v\in L_{V}(S)$ and $w\in S$ and
let $l$ be a nonnegative integer such that
\begin{eqnarray}
(x_{0}+x_{2})^{l}Y(Y(u,x_{0})v,x_{2})w&=&
(x_{0}+x_{2})^{l}Y(u,x_{0}+x_{2})Y(v,x_{2})w\label{e2.62}\\
(x_{0}+x_{2})^{l}Y(u,x_{0})Y(w,-x_{2})v
&=&(x_{0}+x_{2})^{l}Y(w,-x_{2})Y(u,x_{0})v.
\end{eqnarray}
With $v\in L_{V}(S),\; w\in S$, we also have
\begin{eqnarray}\label{e2.64}
Y(v,x_{2})w=e^{x_{2}D}Y(w,-x_{2})v.
\end{eqnarray}
Using (\ref{e2.62})-(\ref{e2.64}) and the conjugation formula 
(\ref{econjugationformula1}) we get
\begin{eqnarray}
(x_{0}+x_{2})^{l}Y(Y(u,x_{0})v,x_{2})w&=&
(x_{0}+x_{2})^{l}Y(u,x_{0}+x_{2})Y(v,x_{2})w\nonumber\\
&=&(x_{0}+x_{2})^{l}Y(u,x_{0}+x_{2})e^{x_{2}D}Y(w,-x_{2})v\nonumber\\
&=&e^{x_{2}D}(x_{0}+x_{2})^{l}Y(u,x_{0})Y(w,-x_{2})v\nonumber\\
&=&e^{x_{2}D}(x_{0}+x_{2})^{l}Y(w,-x_{2})Y(u,x_{0})v.
\end{eqnarray}
Multiplying both sides by $(x_{2}+x_{0})^{-l}$ we obtain
\begin{eqnarray}
Y(Y(u,x_{0})v,x_{2})w=e^{x_{2}D}Y(w,-x_{2})Y(u,x_{0})v,
\end{eqnarray}
as desired. $\;\;\;\;\Box$

A subset $S$ of a weak axiomatic $G_{1}$-vertex algebra $V$ is said to
be {\em local} if $S\subset L_{V}(S)$, that is, 
for any $u,v\in S$, there exists a nonnegative integer $k$ such that
\begin{eqnarray}
(x_{1}-x_{2})^{k}[Y(u,x_{1}),Y(v,x_{2})]=0.
\end{eqnarray}
In view of Remark \ref{rjacobiidentity2}, any local subalgebra of $V$ is an 
ordinary vertex algebra. Furthermore, ordinary vertex algebras
are exactly local weak axiomatic $G_{1}$-vertex algebras.

\bl{lmaximallocalspace}
Any maximal local subspace $A$ of a weak axiomatic $G_{1}$-vertex
algebra $V$ is an ordinary vertex algebra.
\el

{\bf Proof.} It suffices to prove that $A$ is a subalgebra.
Clearly, $A+\C {\bf 1}$ is a local subspace of $V$.
With $A$ being maximal, we must have
$A=A+\C {\bf 1}$, hence ${\bf 1}\in A$.
Let $a,b\in A,\; n\in \Z$. Since $A$ is local, we have $A\subset L_{V}(A)$.
In view of Proposition \ref{plocalizer}, we get $a_{n}b\in L_{V}(A)$.
This also implies that $a,b\in A\subset L_{V}(\{a_{n}b\})$.
By Proposition \ref{plocalizer} again
we have $a_{n}b\in L_{V}(\{a_{n}b\})$. This shows that $A+\C a_{n}b$ is local.
Again, since $A$ is maximal, we must have $a_{n}b\in A$.
This proves that $A$ is a subalgebra of $V$, as we need.
$\;\;\;\;\Box$

\bp{pcommutativesub}
Let $V$ be a weak axiomatic $G_{1}$-vertex algebra and 
let $S$ be a local subset of
$V$. Then the subalgebra $\<S\>$ of $V$ generated by $S$ is 
an ordinary vertex algebra.
\ep

{\bf Proof.} It follows from Zorn's lemma that $S$ is contained in
some maximal local subspace $A$ of $V$. 
By Lemma \ref{lmaximallocalspace},
$A$ is an ordinary vertex algebra, so is $\<S\>$ as a subalgebra of $A$.
$\;\;\;\;\Box$

The following Proposition follows immediately from 
the proof of the corresponding 
proposition in [FHL] for ordinary vertex algebras:

\bp{pproduct}
Let $V_{1},\dots, V_{n}$ be (weak) axiomatic $G_{1}$-vertex algebras. Then 
the tensor product space
$$V_{1}\otimes\cdots \otimes V_{n},$$
equipped with the vertex operator map $Y$ defined by
\begin{eqnarray}
Y(v^{1}\otimes\cdots\otimes v^{n},x)=Y(v^{1},x)\otimes\cdots \otimes Y(v^{n},x)
\end{eqnarray}
and equipped with the vacuum vector
\begin{eqnarray}
{\bf 1}={\bf 1}\otimes\cdots\otimes {\bf 1},
\end{eqnarray}
is a (weak) axiomatic $G_{1}$-vertex algebra. $\;\;\;\;\Box$
\ep

\br{rdefinitionleftrightideal}
{\em The notions of left ideal and right ideal for a weak axiomatic
$G_{1}$-vertex algebra are defined 
in the obvious ways.}
\er

\section{Constructing $G_{1}$-vertex algebras from ordinary vertex algebras}
In this section we shall give several ways to construct 
(weak) axiomatic $G_{1}$-vertex algebras from ordinary vertex algebras.
Specifically, we consider a certain twisting 
of abelian group-graded vertex algebras
by a normalized  $2$-cocycle
and the cross product of a vertex algebra $V$ with a group $G$
which acts on $V$ by automorphisms. 
We also consider the axiomatic $G_{1}$-vertex algebra
$M(n,V)$ of all $n\times n$ matrices over a vertex algebra $V$.
We introduce the notions of
restricted weak axiomatic $G_{1}$-vertex algebra
to unify the given examples.

Let $G$ be an abelian  group. The group algebra $\C[G]$ is 
a commutative associative algebra.
Let $\epsilon(\cdot,\cdot)$ be a {\em normalized 2-cocycle} of $G$
in the sense that $\epsilon(\cdot,\cdot)$ is a $\C^{*}$-valued
function on $G\times G$ such that
\begin{eqnarray}\label{eproperty}
& &\epsilon(\alpha,\beta+\gamma)\epsilon(\beta,\gamma)
=\epsilon(\alpha,\beta)\epsilon(\alpha+\beta,\gamma),\\
& &\epsilon(\alpha,0)=\epsilon(0,\alpha)=1.
\end{eqnarray}
Then $\C[G]$ becomes a noncommutative associative algebra 
(with the same identity element) by defining
\begin{eqnarray}
e^{\alpha}\circ e^{\beta}=\epsilon(\alpha,\beta)e^{\alpha+\beta}
\;\;\;\mbox{ for }\alpha,\beta\in G.
\end{eqnarray}
More generally, let $A$ be a $G$-graded associative algebra. Then,
in the same way we can make $A$ an associative algebra by
using a normalized 2-cocycle $\epsilon$ of $G$.
We next discuss an exact analogue for axiomatic $G_{1}$-vertex algebras of 
this fact. 

\bex{exatwisted}
{\em Let $G$ be an abelian group and let $V$ be a $G$-graded 
axiomatic $G_{1}$-vertex algebra is the sense that 
$V$ is an axiomatic $G_{1}$-vertex algebra equipped with a $G$-grading
$V=\oplus_{g\in G}V^{g}$ such that 
\begin{eqnarray}
u_{n}v\in V^{g+h}\;\;\;\mbox{ for }u\in V^{g},\; v\in V^{h},\; n\in \Z.
\end{eqnarray}
Assume that ${\bf 1}\in V^{0}$. 
Let $\epsilon$ be a normalized $2$-cocycle of $G$. Define a linear map
\begin{eqnarray}
Y_{\epsilon}(\cdot,x): V&\rightarrow& (\End V)[[x,x^{-1}]]\nonumber\\
v &\mapsto& Y_{\epsilon}(v,x)
\end{eqnarray}
by
\begin{eqnarray}
Y_{\epsilon}(u,x)v=\epsilon(g,h)Y(u,x)v\;\;\;\mbox{ for }
u\in V^{g},\; v\in V^{h},\; g,h\in G.
\end{eqnarray}
With $\epsilon$ being normalized, 
it is clear that all the axioms except for the weak associativity hold.
Let $u\in V^{g_{1}},\; v\in V^{g_{2}},\; w\in V^{g_{3}}$. Then
\begin{eqnarray}\label{eproductYepsilon}
Y_{\epsilon}(u,x_{1})Y_{\epsilon}(v,x_{2})w
=\epsilon(g_{1},g_{2}+g_{3})\epsilon(g_{2},g_{3})
Y(u,x_{1})Y(v,x_{2})w
\end{eqnarray}
and
\begin{eqnarray}\label{eiterateYepsilon}
Y_{\epsilon}(Y_{\epsilon}(u,x_{0})v,x_{2})w
=\epsilon(g_{1},g_{2})\epsilon(g_{1}+g_{2},g_{3})Y(Y(u,x_{0})v,x_{2})w.
\end{eqnarray}
With the property (\ref{eproperty}), we easily see that 
the weak associativity holds for $Y_{\epsilon}$.
Thus, $(V,Y_{\epsilon},{\bf 1})$ is an axiomatic $G_{1}$-vertex algebra.

Furthermore, we assume that $V$ is an ordinary vertex algebra.
{}From (\ref{eproductYepsilon}) (using the
obvious symmetry) we have
\begin{eqnarray}\label{eproductYepsilon2}
Y_{\epsilon}(v,x_{2})Y_{\epsilon}(u,x_{1})w
=\epsilon(g_{2},g_{1}+g_{3})\epsilon(g_{1},g_{3})
Y(v,x_{2})Y(u,x_{1})w.
\end{eqnarray}
Then it follows immediately from
(\ref{eproductYepsilon})-(\ref{eproductYepsilon2}),
(\ref{eproperty}) and the Jacobi identity (\ref{ejacobidef}) that 
\begin{eqnarray}\label{erestrictedjacobi1}
& &x_{0}^{-1}
\delta\left(\frac{x_{1}-x_{2}}{x_{0}}\right)Y(u,x_{1})Y(v,x_{2})
-c(g,h)x_{0}^{-1}
\delta\left(\frac{x_{2}-x_{1}}{-x_{0}}\right)Y(v,x_{2})Y(u,x_{1})\nonumber\\
&=&x_{2}^{-1}\delta\left(\frac{x_{1}-x_{0}}{x_{2}}\right)Y(Y(u,x_{0})v,x_{2})
\end{eqnarray}
for $u\in V^{g},\; v\in V^{h},\; g,h\in G$, where
\begin{eqnarray}
c(g,h)=\epsilon(g,h)\epsilon(h,g)^{-1}.
\end{eqnarray}
It is easy to see that $c(\cdot,\cdot)$ is a $\C^{*}$-valued 
bilinear form on $G$.}
\eex

\br{rrestrictedgva1}
{\em Note that for an ordinary vertex algebra $V$, 
the axiomatic $G_{1}$-vertex algebra $(V,Y_{\epsilon}, {\bf 1})$ 
obtained in Example \ref{exatwisted} in fact has a natural generalized 
vertex algebra structure
$(V,Y_{\epsilon}, {\bf 1}, G, c, (\cdot,\cdot))$ with $(\cdot,\cdot)=0$
in the sense of [DL]. The notion of generalized vertex algebra
naturally generalizes the notions of ordinary vertex algebra and
vertex superalgebra. In the notion of
generalized vertex algebra $(V,Y, {\bf 1}, G, c, (\cdot,\cdot))$
defined in [DL], if the form $(\cdot,\cdot)$ is (zero) trivial, then
the linear map $Y$ maps $V$ to $\Hom (V,V((x)))$ and the generalized
Jacobi identity for the generalized vertex algebra reduces to
(\ref{erestrictedjacobi1}). For convenience,
we call a generalized vertex algebra with trivial form $(\cdot,\cdot)$
a {\em restricted generalized vertex algebra}.
It is easy to see that any restricted generalized vertex algebra is an
axiomatic $G_{1}$-vertex algebra.}
\er

%Lie color algebra [X]?

\br{rexadllattice}
{\em Let $L$ be an integral lattice. For any 
normalized $2$-cocycle $\epsilon(\cdot,\cdot)$ of $L$,
a restricted generalized vertex algebra $V_{L}$ was constructed in [DL].
It was proved in [DL] that for a certain $\epsilon(\cdot,\cdot)$,
$V_{L}$ is a vertex superalgebra.}
\er

%\bex{exaaffine}
%{\em Let ${\fg}$ be a finite-dimensional simple Lie algebra over $\C$ 
%equipped with the normalized Killing form $\<\cdot,\cdot\>$.
%Denote the affine algebra $\hat{\fg}$ (without $d$).
%For any complex number $\ell$, let $\C_{\ell}$ be 
%the $1$-dimensional $({\fg}\otimes \C[t]+\C c)$-module 
%with ${\fg}\otimes \C[t]$ acting as zero
%and with $c$ acting as scalar $\ell$. Form the induced $\hat{\fg}$-module
%\begin{eqnarray}
%M_{\fg}(\ell,0)=U(\hat{\fg})\otimes _{U({\fg}\otimes \C[t]+\C c)}\C_{\ell}.
%\end{eqnarray}
%Then $M_{\fg}(\ell,0)$ has a natural vertex (operator) algebra
%structure ([FF], [FZ]).

%Let ${\bf H}$ be a Cartan subalgebra of ${\fg}$ and 
%let ${\bf Q}$ be the root lattice of ${\fg}$ with respect to ${\bf H}$.
%Then $M_{\fg}(\ell,0)$ is naturally ${\bf Q}$-graded.
%For any normalized $2$-cocycle $\epsilon$ on ${\bf Q}$, we have 
%a $G_{1}$-vertex 
%algebra structure $Y_{\epsilon}(\cdot,x)$ on $M_{\fg}(\ell,0)$.
%Denote by $L_{\fg}(\ell,0)$ the irreducible quotient module 
%of $M_{\fg}(\ell,0)$. Then $L_{\fg}(\ell,0)$ is 
%a vertex (operator) algebra and it is also ${\bf Q}$-graded.}
%\eex

\bex{examplematrix}
{\em In view of Proposition \ref{pproduct} and 
Example \ref{exampleassociativealgebra},
the tensor product of any axiomatic $G_{1}$-vertex algebra $V$ 
with any associative algebra $A$ is an axiomatic $G_{1}$-vertex algebra. 
In particular, the tensor product $V\otimes M(n,\C)$ is an 
axiomatic $G_{1}$-vertex algebra. 
If we naturally identify the space $V\otimes M(n,\C)$
with the vector space $M(n,V)$ of $n\times n$ matrices with entries in $V$,
then for $A=(a_{ij}),\; B\in M(n,V)$,
\begin{eqnarray}
Y(A,x)B=(Y(a_{ij},x))B \;\;\;\mbox{(the formal matrix product)}.
\end{eqnarray}
Furthermore, the general linear group $GL(n,\C)$ naturally acts 
on $M(n,V)$, and so does the unitary group $U(n)$.}
\eex

In the classical associative algebra theory, for an algebra $A$ 
acted by a group $G$,
there is a notion of cross product (or skew group algebra)
where the underlying vector space is
$A\otimes \C[G]$ and the multiplication is given by
\begin{eqnarray}
(ag)(bh)=ag(b)gh\;\;\;\mbox{ for }a,b\in A,\; g,h\in G.
\end{eqnarray}
When $G$ acts on $A$ trivially,  the cross product algebra becomes 
the usual product of $A$ with the group algebra $\C[G]$.
In the following we give an analogue for axiomatic $G_{1}$-vertex algebras.
To do so, first we define the notion of automorphism of a (weak)
axiomatic $G_{1}$-vertex algebra 
in the obvious way: An {\em automorphism} of $V$ is an invertible linear 
endomorphism $\psi$ of $V$ such that $\psi({\bf 1})={\bf 1}$ and such that
$\psi(u_{n}v)=\psi(u)_{n}\psi(v)$ for $u,v\in V,\; n\in \Z$.

\bex{examplegroupalgebracrossproduct}
{\em Let $V$ be an axiomatic $G_{1}$-vertex algebra and 
let $G$ be a group acting on $V$ as automorphisms. 
Define a vertex operator map $Y$ on $V[G]=V\otimes \C[G]$ by
\begin{eqnarray}
Y(ug,x)(vh)=Y(u,x)g(v)gh\;\;\;\mbox{ for }u,v\in V,\; g,h\in G.
\end{eqnarray}
Taking ${\bf 1}e$ to be the vacuum vector, where $e$ denotes 
the identity element of $G$, 
we easily see that all the axioms except for the weak associativity hold.
Furthermore, let $u,v,w\in V,\; g_{1},g_{2},g_{3}\in G$. 
Since $G$ acts as automorphisms of $V$, we have
\begin{eqnarray}\label{ecrossproduct1}
Y(ug_{1},x_{1})Y(vg_{2},x_{2})wg_{3}
&=&Y(ug_{1},x_{1})Y(v,x_{2})g_{2}(w)(g_{2}g_{3})\nonumber\\
&=&Y(u,x_{1})g_{1}(Y(v,x_{2})g_{2}(w))g_{1}(g_{2}g_{3})\nonumber\\
&=&Y(u,x_{1})Y(g_{1}(v),x_{2})(g_{1}g_{2})(w)g_{1}(g_{2}g_{3})
\end{eqnarray}
and
\begin{eqnarray}
Y(Y(ug_{1},x_{0})vg_{2},x_{2})wg_{3}
&=&Y(Y(u,x_{0})g_{1}(v)g_{1}g_{2},x_{2})wg_{3}\nonumber\\
&=&Y(Y(u,x_{0})g_{1}(v),x_{2})(g_{1}g_{2})(w)(g_{1}g_{2})g_{3}.
\end{eqnarray}
Now, the weaker version of the weak associativity for $Y$ on $V[G]$ 
follows immediately.
Therefore,  $V[G]$ equipped with this vertex operator map $Y$
is a weak axiomatic $G_{1}$-vertex algebra. 
%We shall denote this $G_{1}$-vertex algebra by $V\rtimes G$. 
It is clear that if $\dim \C[G]w <\infty$ for
every $w\in V$, $V[G]$ is an axiomatic $G_{1}$-vertex algebra. 
In particular, this is true if $V$ is an ordinary vertex
operator algebra (with the two grading restrictions), since 
any automorphism group preserves every homogeneous subspace of $V$.}
\eex

\bex{examplegroupalgebracrossproductcontinuation}
{\em We continue with Example \ref{examplegroupalgebracrossproduct}
studying the cross product weak axiomatic $G_{1}$-vertex algebra.
Let us assume that $V$ is an ordinary vertex algebra.
{}From (\ref{ecrossproduct1}) (using the obvious symmetry) we have
\begin{eqnarray}
Y(vg_{2},x_{2})Y(ug_{1},x_{1})wg_{3}
=Y(v,x_{2})Y(g_{2}(u),x_{1})(g_{2}g_{1})(w) g_{2}(g_{1}g_{3}).
\end{eqnarray}
We are not able to see weak commutativity.
On the other hand, notice that
\begin{eqnarray}
& &Y(g_{1}(v)g_{2},x_{2})Y(g_{2}^{-1}(u)g_{1},x_{1})
(g_{1}^{-1}g_{2}^{-1}g_{1}g_{2})(w)
(g_{1}^{-1}g_{2}^{-1}g_{1}g_{2}g_{3})\nonumber\\
&=&Y(g_{1}(v),x_{2})Y(u,x_{1})(g_{1}g_{2})(w) (g_{1}g_{2}g_{3}).
\end{eqnarray}
Then we have the following Jacobi-like identity:
\begin{eqnarray}
& &x_{0}^{-1}\delta\left(\frac{x_{1}-x_{2}}{x_{0}}\right)
Y(ug_{1},x_{1})Y(vg_{2},x_{2})wg_{3}\nonumber\\
& &-x_{0}^{-1}\delta\left(\frac{x_{2}-x_{1}}{-x_{0}}\right)
Y(g_{1}(v)g_{2},x_{2})Y(g_{2}^{-1}(u)g_{1},x_{1})
(g_{1}^{-1}g_{2}^{-1}g_{1}g_{2})(w)
(g_{1}^{-1}g_{2}^{-1}g_{1}g_{2}g_{3})\nonumber\\
&=&x_{2}^{-1}\delta\left(\frac{x_{1}-x_{0}}{x_{2}}\right)
Y(Y(ug_{1},x_{0})vg_{2},x_{2})wg_{3}.
\end{eqnarray}
Let $R$ be the linear endomorphism of $V[G]\otimes V[G]\otimes V[G]$ 
uniquely determined by
\begin{eqnarray}
R(vg_{2}\otimes ug_{1}\otimes wg_{3})
=g_{1}(v)g_{2}\otimes g_{2}^{-1}(u)g_{1}\otimes 
(g_{1}^{-1}g_{2}^{-1}g_{1}g_{2})(w)
(g_{1}^{-1}g_{2}^{-1}g_{1}g_{2}g_{3}).
\end{eqnarray}
Then 
\begin{eqnarray}
& &x_{0}^{-1}\delta\left(\frac{x_{1}-x_{2}}{x_{0}}\right)
Y(ug_{1},x_{1})Y(vg_{2},x_{2})wg_{3}\nonumber\\
& &-x_{0}^{-1}\delta\left(\frac{x_{2}-x_{1}}{-x_{0}}\right)
(Y\otimes Y)(x_{2},x_{1})R(vg_{2}\otimes ug_{1}\otimes wg_{3})
\nonumber\\
&=&x_{2}^{-1}\delta\left(\frac{x_{1}-x_{0}}{x_{2}}\right)
Y(Y(ug_{1},x_{0})vg_{2},x_{2})wg_{3},
\end{eqnarray}
where $(Y\otimes Y)(x_{1},x_{2})$ is the linear map from 
$V[G]\otimes V[G]\otimes V[G]$ to $V[G]((x_{1}))((x_{2}))$ defined by
\begin{eqnarray}
(Y\otimes Y)(x_{1},x_{2})(a\otimes b\otimes c)=Y(a,x_{1})Y(b,x_{2})c.
\end{eqnarray}
If $G$ is abelian, the Jacobi-like identity reduces to
\begin{eqnarray}
& &x_{0}^{-1}\delta\left(\frac{x_{1}-x_{2}}{x_{0}}\right)
Y(ug_{1},x_{1})Y(vg_{2},x_{2})wg_{3}\nonumber\\
& &-x_{0}^{-1}\delta\left(\frac{x_{2}-x_{1}}{-x_{0}}\right)
Y(g_{1}(v)g_{2},x_{2})Y(g_{2}^{-1}(u)g_{1},x_{1})
wg_{3}\nonumber\\
&=&x_{2}^{-1}\delta\left(\frac{x_{1}-x_{0}}{x_{2}}\right)
Y(Y(ug_{1},x_{0})vg_{2},x_{2})wg_{3}.
\end{eqnarray}}
\eex

\br{rtwistedcrossproduct}
{\em The classical cross product algebra or skew group algebra
was employed by Yamskulna [Y] in a certain study of twisted modules. 
It would be interesting if one can apply the cross product 
axiomatic $G_{1}$-vertex algebra in that type of study.}
\er

Motivated by Example \ref{examplegroupalgebracrossproductcontinuation},
we introduce the following notion of 
restricted weak axiomatic $G_{1}$-vertex algebra:

\bd{dquasinoncommva}
{\em A {\em restricted weak axiomatic $G_{1}$-vertex algebra} 
is a vector space $V$ 
equipped with a linear map $Y$ from $V$ to $\Hom (V,V((x)))$,
a distinguished vector ${\bf 1}\in V$ and a linear map 
$R\in \End (V\otimes V\otimes V)$ such that 
$Y({\bf 1},x)=1$, for $v\in V$,
$Y(v,x){\bf 1}\in V[[x]]$ and $\lim_{x\rightarrow 0}Y(v,x){\bf 1}=v$
and such that the following {\em Jacobi-like identity} holds 
for $u,v,w\in V$:
\begin{eqnarray}\label{equasijacobi}
& &x_{0}^{-1}
\delta\left(\frac{x_{1}-x_{2}}{x_{0}}\right)Y(u,x_{1})Y(v,x_{2})w\nonumber\\
& &-x_{0}^{-1}
\delta\left(\frac{x_{2}-x_{1}}{-x_{0}}\right)
(Y\otimes Y)(x_{2},x_{1})R(v\otimes u\otimes w)
\nonumber\\
&=&x_{2}^{-1}\delta\left(\frac{x_{1}-x_{0}}{x_{2}}\right)
Y(Y(u,x_{0})v,x_{2})w,
\end{eqnarray}
where $(Y\otimes Y)(x_{1},x_{2})$ is the linear map from 
$V\otimes V\otimes V$ to $V((x_{1}))((x_{2}))$ defined by
\begin{eqnarray}
(Y\otimes Y)(x_{1},x_{2})(u,v,w)=Y(u,x_{1})Y(v,x_{2})w.
\end{eqnarray}}
\ed

\br{rquasi}
{\em Let $u,v,w\in V$ and let $l$ be a nonnegative integer such that
$$x_{1}^{l}(Y\otimes Y)(x_{2},x_{1})R(v\otimes u\otimes w)
\in V[[x_{1},x_{2},x_{2}^{-1}]].$$
By taking $\Res_{x_{1}}x_{1}^{l}$ of the
Jacobi-like identity (\ref{equasijacobi}) we get
\begin{eqnarray}
(x_{0}+x_{2})^{l}Y(u,x_{0}+x_{2})Y(v,x_{2})w
=(x_{0}+x_{2})^{l}Y(Y(u,x_{0})v,x_{2})w.
\end{eqnarray}
Then, a restricted weak axiomatic $G_{1}$-vertex algebra 
is indeed a weak axiomatic $G_{1}$-vertex algebra.
Furthermore, let $k$ be a nonnegative integer such that
$x^{k}Y(u,x)v\in V[[x]]$. Then we get
\begin{eqnarray}
(x_{1}-x_{2})^{k}Y(u,x_{1})Y(v,x_{2})w=(x_{1}-x_{2})^{k}
(Y\otimes Y)(x_{2},x_{1})R(v\otimes u\otimes w).
\end{eqnarray}}
\er

\br{rquasimany}
{\em We explain that many examples of weak axiomatic $G_{1}$-vertex
algebras, discussed previously, are restricted
weak axiomatic $G_{1}$-vertex algebras.
First, from Example \ref{examplegroupalgebracrossproductcontinuation}
the cross product of an ordinary vertex
algebra with a group is a restricted weak axiomatic $G_{1}$-vertex algebra.

Second, any restricted generalized vertex algebra 
defined in Remark \ref{rrestrictedgva1} is a 
restricted weak axiomatic $G_{1}$-vertex algebra, where
\begin{eqnarray}
R(v\otimes u\otimes w)=c(g,h)(v\otimes u\otimes w)
\;\;\;\mbox{ for }u\in V^{g},\;v\in V^{h},\;w\in V,\;g,h\in G.
\end{eqnarray}

Third, for any ordinary vertex algebra $V$ 
and any associative algebra $A$,
the tensor product axiomatic $G_{1}$-vertex algebra $V\otimes A$ 
is a restricted weak axiomatic $G_{1}$-vertex algebra, where
\begin{eqnarray}
R(ua\otimes vb\otimes wc)=ub \otimes va\otimes wc
\;\;\;\mbox{ for }u,v,w\in V,\; a,b,c\in A.
\end{eqnarray}
Indeed, the Jacobi-like identity holds because of
the Jacobi identity of $V$ and the relations
\begin{eqnarray}
& &Y(ua,x_{1})Y(vb,x_{2})wc=Y(u,x_{1})Y(v,x_{2})wa(bc),\\
& &(Y\otimes Y)(x_{2},x_{1})R(vb\otimes ua\otimes wc)=
Y(va,x_{2})Y(ub,x_{1})wc=Y(v,x_{2})Y(u,x_{1})wa(bc),\;\;\;\;\;\\
& &Y(Y(ua,x_{0})vb,x_{2})wc=Y(Y(u,x_{0})v,x_{2})w(ab)c.
\end{eqnarray}
Then, for any ordinary vertex algebra $V$, $M(n,V)$ is a 
restricted weak axiomatic $G_{1}$-vertex algebra.
On the other hand, by taking $V=\C$, we see that any 
associative algebra $A$ is a restricted weak axiomatic 
$G_{1}$-vertex algebra,
where the classical associativity axiom is expressed 
in terms of delta functions as
\begin{eqnarray}
x_{0}^{-1}
\delta\left(\frac{x_{1}-x_{2}}{x_{0}}\right)a(bc)
-x_{0}^{-1}
\delta\left(\frac{x_{2}-x_{1}}{-x_{0}}\right)a(bc)
=x_{2}^{-1}\delta\left(\frac{x_{1}-x_{0}}{x_{2}}\right)(ab)c.
\end{eqnarray}}
\er

\br{rslocalva}
{\em Since restricted weak axiomatic 
$G_{1}$-vertex algebras are still too general to study,
we may consider a special class of restricted weak axiomatic 
$G_{1}$-vertex algebras with $R=R'\otimes 1$, where $R'\in \End
(V\otimes V)$. Restricted generalized vertex algebras, matrix
axiomatic $G_{1}$-vertex algebras $M(n,V)$ and the cross
product axiomatic $G_{1}$-vertex algebras $V[G]$ with $G$ abelian are
of this type.}
\er

%\bd{dquasi}
%{\em A {\em quasi-local axiomatic $G_{1}$-vertex algebra} 
%is a restricted weak axiomatic $G_{1}$-vertex algebra $V$
%for which the linear map $R\in \End (V\otimes V\otimes V)$ is given by
%\begin{eqnarray}
%R(u\otimes v\otimes w)
%=\sum_{\alpha\in I}f_{\alpha}(v)\otimes g_{\alpha}(u)\otimes w
%\end{eqnarray}
%for $u,v,w\in V$,
%where $I$ is an index set and $f_{\alpha},g_{\alpha}\in \End V$ such that
%for every $u\in V$, $g_{\alpha}(u)\ne 0$ for only finitely many $\alpha$.}
%\ed

%\br{rquasilocalva}
%{\em Let $V$ be a quasi-local axiomatic $G_{1}$-vertex algebra.
%For any $u,w\in V$, there exists a nonnegative integer $l$ such that
%\begin{eqnarray}
%x^{l}Y(g_{\alpha}(u),x)w\in V[[x]]\;\;\;\mbox{ for all }\alpha\in I.
%\end{eqnarray}
%Then from the Jacobi-like identity we get
%\begin{eqnarray}
%(x_{0}+x_{2})^{l}Y(u,x_{0}+x_{2})Y(v,x_{2})w
%=(x_{0}+x_{2})^{l}Y(Y(u,x_{0})v,x_{2})w
%\end{eqnarray}
%for all $v\in V$. Thus $V$ is an axiomatic $G_{1}$-vertex algebra.
%On the other hand, for any $u,v\in V$, let $k$ be a nonnegative integer
%such that $x^{k}Y(u,x)v\in V[[x]]$. Then from the Jacobi-like identity
%we get
%\begin{eqnarray}
%(x_{1}-x_{2})^{k}Y(u,x_{1})Y(v,x_{2})
%=\sum_{\alpha\in I}(x_{1}-x_{2})^{k}
%Y(f_{\alpha}(v),x_{2})Y(g_{\alpha}(u),x_{1}),
%\end{eqnarray}
%which is a generalization of the weak commutativity relation, or
%locality relation.
%This justifies the designation of the notion of quasi-local 
%axiomatic $G_{1}$-vertex algebra. 
%For modules, this relation also holds.}
%\er

\br{rag1voadefinition}
{\rm We here define a notion of axiomatic $G_{1}$-vertex operator algebra.
An {\em axiomatic $G_{1}$-vertex operator algebra} is 
an axiomatic $G_{1}$-vertex algebra $V$
equipped with a distinguished vector $\omega\in V$, called the 
{\em conformal vector}, such that
\begin{eqnarray}
[L(m),L(n)]=(m-n)L(m+n)+{1\over 12}(m^{3}-m){\rm rank}V \delta_{m+n,0}
\end{eqnarray}
for $m,n\in \Z$, 
where
\begin{eqnarray}
Y(\omega,x)=\sum_{n\in \Z}L(n)x^{-n-2}
\end{eqnarray}
and ${\rm rank}V\in \C$ is called the {\em rank} of $V$, and such that
\begin{eqnarray}
& &Y(L(-1)v,x)={d\over dx}Y(v,x),\\
& &[L(m),Y(v,x)]=\sum_{i\ge 0}{m+1\choose i}x^{m+1-i}Y(L(i-1)v,x)
\end{eqnarray}
for $v\in V$ and such that
\begin{eqnarray}
V=\coprod_{n\in \Z}V_{(n)},
\end{eqnarray}
$\dim V_{(n)}<\infty$ for all $n\in \Z$ and
$V_{(n)}=0$ for $n$ sufficiently negative,
where for $n\in \Z$, 
\begin{eqnarray}
V_{(n)}=\{ v\in V\;|\; L(0)v=nv\}.
\end{eqnarray}}
\er

%Equivalently, for every $v\in V$, there exists a nonnegative 
%integer $k$ such that
%\begin{eqnarray}
%(x_{1}-x_{2})^{k}[ Y(\omega,x_{1}),Y(v,x_{2})]=0.
%\end{eqnarray}

%Note that we do not have examples of axiomatic 
%$G_{1}$-vertex operator algebras other than those from
%the three constructions.
%Construct essentially new axiomatic $G_{1}$-vertex operator algebras.
%Simple and $\ker D=\C {\bf 1}$.
%More generally, change Heisenberg algebra.
%boson, fermion, quon operators.

%The theory of contragredient module established in [FHL] 
%for vertex operator algebras can be extended to 
%conformal $G_{1}$-vertex operator algebras ?

%Let $V$ be a conformal $G_{1}$-vertex operator algebra 
%and let $W=\oplus_{h\in \C}W_{(h)}$ be a $V$-module.
%Set 
%\begin{eqnarray}
%W'=\oplus_{h\in \C}W_{(h)}^{*}
%\end{eqnarray}
%and define $Y'(v,x)$ as usual:
%\begin{eqnarray}
%\<Y'(v,x)w',w\>=\<w',Y(e^{xL(1)}(-x^{-2})^{L(0)}v,x^{-1})\>
%\end{eqnarray}
%for $w'\in W,\; w\in W$. Then
%Let $u,v\in V$, $w'\in W',\; w\in W$, we have
%\begin{eqnarray}
%& &\<Y'(u,x_{0}+x_{2})Y'(v,x_{2})w',w\>\nonumber\\
%&=&\<w', Y(e^{x_{2}L(1)}(-x_{2}^{-2})^{L(0)}v,x_{2}^{-1})
%Y(e^{(x_{0}+x_{2})L(1)}(-(x_{0}+x_{2})^{-2})^{L(0)}u,(x_{0}+x_{2})^{-1})w\>.
%\end{eqnarray}
%(it is not true for an associative algebra.)

\section{Modules for axiomatic $G_{1}$-vertex algebras}
In this section we shall define the notion of module
for a (weak) axiomatic $G_{1}$-vertex algebra and we obtain
certain analogous results of ordinary vertex algebras
for axiomatic $G_{1}$-vertex algebras.

Let $V$ be a weak axiomatic $G_{1}$-vertex algebra, fixed throughout
this section. We first define the notion of $V$-module.

\bd{dmodule}
{\em A $V$-{\em module} is a vector space $W$ equipped with a linear map
\begin{eqnarray}
Y_{W}(\cdot,x): V &\rightarrow & (\End W)[[x,x^{-1}]]\nonumber\\
v&\mapsto& Y_{W}(v,x)=\sum_{n\in \Z}v_{n}x^{-n-1}\;\;\; (v_{n}\in \End W)
\end{eqnarray}
such that all the following axioms hold: For every $v\in V,\; w\in W$, 
\begin{eqnarray}
v_{n}w=0\;\;\;\mbox{ for $n$ sufficiently large};
\end{eqnarray}
\begin{eqnarray}
Y_{W}({\bf 1},x)=1_{W}\;\;\;
\mbox{(where $1_{W}$ is the identity operator on $W$)};
\end{eqnarray}
for any $u, v\in V$ and $w\in W$, there exists $l\in \N$
such that
\begin{eqnarray}\label{emoduleweakassoc}
(x_{0}+x_{2})^{l}Y_{W}(u,x_{0}+x_{2})Y_{W}(v,x_{2})w
=(x_{0}+x_{2})^{l}Y_{W}(Y(u,x_{0})v,x_{2})w.
\end{eqnarray}}
\ed

If $V$ is an axiomatic $G_{1}$-vertex algebra,
for the notion of $V$-module, we use the stronger weak associativity:
For any $u\in V$ and $w\in W$, there exists $l\in \N$
such that for {\em all} $v\in V$, (\ref{emoduleweakassoc}) holds

We next discuss  some consequences of the definition.
First, by carefully examining the first half of the proof of 
Proposition \ref{pcon1} we find that 
the same argument (with $v$ being replaced by $w\in W$) gives:

\bp{pmoduleDproperty}
Let $(W,Y_{W})$ be a $V$-module. Then
\begin{eqnarray}
Y_{W}(Dv,x)={d\over dx}Y_{W}(v,x)\;\;\;\mbox{ for }v\in V
\end{eqnarray}
(recall the linear operator $D$ on $V$).$\;\;\;\;\Box$
\ep

The notions of submodule, irreducible module and module homomorphism
are defined in the obvious ways.

The following result tells us how a certain commutativity
relation of vertex operators on $V$ is related to that of the
vertex operators on other modules:

\bp{pcommutant}
Let $u,v\in V,\; k\in \N$ and $q\in \C^{*}$. If
\begin{eqnarray}\label{eqlocalityonalgebra}
(x_{1}-x_{2})^{k}Y(u,x_{1})Y(v,x_{2})=q
(x_{1}-x_{2})^{k}Y(v,x_{2})Y(u,x_{1})
\end{eqnarray}
on $V$, then for any $V$-module $(W,Y_{W})$,
\begin{eqnarray}\label{eqlocalityonmodule}
(x_{1}-x_{2})^{k}Y_{W}(u,x_{1})Y_{W}(v,x_{2})=q
(x_{1}-x_{2})^{k}Y_{W}(v,x_{2})Y_{W}(u,x_{1}).
\end{eqnarray}
In particular, if $[Y(u,x_{1}),Y(v,x_{2})]=0$ on $V$, 
then $[Y_{W}(u,x_{1}),Y_{W}(v,x_{2})]=0$ on $W$.
On the other hand, if $W$ is a faithful module and if 
(\ref{eqlocalityonmodule}) holds,
then (\ref{eqlocalityonalgebra}) holds.
\ep

{\bf Proof.} From Proposition \ref{pcommutativity} we have
\begin{eqnarray}\label{etruncation-skewsymmetry}
x^{k}Y(u,x)v\in V[[x]]\;\;\mbox{  and }Y(u,x)v=qe^{xD}Y(v,-x)u.
\end{eqnarray}
Then (\ref{eqlocalityonmodule}) follows from the same proof of 
the ``if'' part of Proposition \ref{pcommutativity}. 
On the other hand, if $W$ is a faithful module and if 
(\ref{eqlocalityonmodule}) holds, the same proof of the
``only if'' part of Proposition \ref{pcommutativity}
shows that (\ref{etruncation-skewsymmetry}) holds.
Then (\ref{eqlocalityonalgebra}) follows from (the ``if'' part of)
Proposition \ref{pcommutativity}.
$\;\;\;\;\Box$

\bp{prestrictedgvamodule}
Let $V$ be a restricted generalized vertex algebra and let $(W,Y_{W})$ 
be a module for $V$ viewed as an axiomatic $G_{1}$-vertex algebra. Then
\begin{eqnarray}
& &x_{0}^{-1}\delta\left(\frac{x_{1}-x_{2}}{x_{0}}\right)
Y_{W}(u,x_{1})Y_{W}(v,x_{2})
-c(g,h)x_{0}^{-1}\delta\left(\frac{x_{2}-x_{1}}{-x_{0}}\right)
Y_{W}(v,x_{2})Y_{W}(u,x_{1})
\nonumber\\
&=&x_{2}^{-1}\delta\left(\frac{x_{1}-x_{0}}{x_{2}}\right)
Y_{W}(Y(u,x_{0})v,x_{2})
\end{eqnarray}
for $u\in V^{g},\; v\in V^{h}$ with $g,h\in G$.
\ep

{\bf Proof.} For $u\in V^{g},\; v\in V^{h}$,
by Proposition \ref{pcommutant}, we have
\begin{eqnarray}
(x_{1}-x_{2})^{k}Y_{W}(u,x_{1})Y_{W}(v,x_{2})=c(g,h)
(x_{1}-x_{2})^{k}Y_{W}(v,x_{2})Y_{W}(u,x_{1}).
\end{eqnarray}
Then it follows immediately from Lemma \ref{lformaljacobiidentity}.
$\;\;\;\;\Box$

%For quasi-local $G_{1}$-vertex algebras?

The following Proposition follows from the proof of
the corresponding result in [FHL] for ordinary vertex algebras.

\bp{ptensorproductmodule}
Let $V_{1},\dots, V_{r}$ be (weak) axiomatic $G_{1}$-vertex algebras and
let $W_{i}$ be a $V_{i}$-module for $i=1,\dots, r$. Then
$W_{1}\otimes \cdots\otimes W_{r}$ is a 
$V_{1}\otimes\cdots\otimes V_{r}$-module with the vertex operator map
$Y$ defined by
\begin{eqnarray}
Y(v^{(1)}\otimes\cdots\otimes v^{(r)},x)
=Y(v^{(1)},x)\otimes\cdots \otimes Y(v^{(r)},x). \;\;\;\;\Box
\end{eqnarray}
\ep

\bex{examatrix}
{\em Let $V$ be a (weak) axiomatic $G_{1}$-vertex algebra, $W$ a $V$-module and
let $n$ be a positive integer. Then $W^{n}$ is an $M(n,V)$-module.}
\eex

As an immediate consequence of Proposition \ref{pcommutant} we have:

\bc{ctensormodule}
Let $U$ and $V$ be (weak) axiomatic $G_{1}$-vertex algebras and let $W$ be a
module for the tensor product $U\otimes V$. Then $W$ is a natural
$U$-module and $V$-module and furthermore, the actions of $U$ and $V$ 
on $W$ commute. $\;\;\;\;\Box$
\ec

Now, we have:

\bp{p}
Let $V$ be an axiomatic $G_{1}$-vertex algebra and 
let $n$ be a positive integer. Then
any irreducible $M(n,V)$-module is of the form $W^{n}$, 
where $W$ is an irreducible $V$-module. On the other hand, 
for any irreducible $V$-module $W$, $W^{n}$ is an irreducible
 $M(n,V)$-module.
\ep

{\bf Proof.} Note that  $\C^{n}$ is the only irreducible module
up to equivalence for the matrix algebra $M(n,\C)$ and 
any $M(n,\C)$-module is completely reducible.
Since any $M(n,V)$-module $W$ is naturally an $M(n,\C)$-module,
we have the canonical decomposition
\begin{eqnarray}
M=\Hom_{M(n,\C)}(\C^{n},M)\otimes \C^{n},
\end{eqnarray}
where $\Hom_{M(n,\C)}(\C^{n},M)$ is naturally a $V$-module.
If $M$ is an irreducible $M(n,V)$-module,
$\Hom_{M(n,\C)}(\C^{n},M)$ is necessarily an irreducible $V$-module.

On the other hand, let $W$ be an irreducible $V$-module.
For $1\le i,j\le n$, denote by $E_{ij}$ be the matrix
whose entry is $1$ at $ij$-position and is zero elsewhere.
Also, for $1\le i\le n$, denote by $e_{i}$ the 
element of $\C^{n}$ whose $i$-th entry is $1$ and others are zero.
Then
$$W^{n}=W\otimes \C^{n}=\sum_{i=1}^{n}W\otimes e_{i}$$
and $E_{ii}W^{n}=W\otimes e_{i}$ for $i=1,\dots,n$.
Since $W$ is an irreducible $V$-module, for $1\le i\le n$,
any nonzero element of $W\otimes e_{i}$ generates $W^{n}$ 
as an $M(n,V)$-module. For any nonzero $w\in W^{n}$,
since $w=E_{11}w+\cdots +E_{nn}w$, $E_{ii}w\ne 0$ for some $i$.
Then it follows that any nonzero element $w$ of $W^{n}$
generates $W^{n}$ as an $M(n,V)$-module.
That is, $W^{n}$ is an irreducible $M(n,V)$-module.
$\;\;\;\;\Box$

%\br{rmarita}
%{\em It can be easily shown that the category of $V$-modules and 
%the category of $M(n,V)$-modules are Morita equivalent.}
%\er

Next, we shall derive a certain compatibility of vertex operators 
$Y_{W}(v,x)$ for $v\in V$ and for a $V$-module $(W,Y_{W})$. 
First we have:

\bl{lcompat2vo}
Let $(W,Y_{W})$ be a $V$-module and
let $u,v\in V$. Then there exists a 
nonnegative integer $k$ such that
\begin{eqnarray}\label{ecompat2vo}
(x_{1}-x_{2})^{k}Y_{W}(u,x_{1})Y_{W}(v,x_{2})\in \Hom (W,W((x_{1},x_{2}))).
\end{eqnarray}
\el

{\bf Proof.} Let $w\in W$. Then there exists a nonnegative integer $l$ 
such that
\begin{eqnarray}\label{eweakassocmoduleconsequence}
(x_{0}+x_{2})^{l}Y_{W}(u,x_{0}+x_{2})Y_{W}(v,x_{2})w=
(x_{2}+x_{0})^{l}Y_{W}(Y(u,x_{0})v,x_{2})w.
\end{eqnarray}
Just as in [FHL], [DL], or [LL], we notice that
the expression on the left-hand side of
(\ref{eweakassocmoduleconsequence}) involves only finitely many
negative powers of $x_{2}$ and the expression on the right-hand side
involves only finitely many
negative powers of $x_{0}$. Consequently, the common quantity
lies in $W((x_{0},x_{2}))$. Let $k$ be a nonnegative integer such that
$x^{k}Y(u,x)v\in V[[x]]$. (Of course, $k$ depends only on $u$ and $v$.) Set
\begin{eqnarray}
p(x_{0},x_{2})=x_{0}^{k}(x_{2}+x_{0})^{l}Y_{W}(Y(u,x_{0})v,x_{2})w.
\end{eqnarray}
Then $p(x_{0},x_{2})\in W[[x_{0},x_{2}]][x_{2}^{-1}]$ and
\begin{eqnarray}
x_{0}^{k}(x_{0}+x_{2})^{l}Y_{W}(u,x_{0}+x_{2})Y_{W}(v,x_{2})w=p(x_{0},x_{2}).
\end{eqnarray}
Applying $e^{-x_{2}{\partial\over \partial x_{0}}}$ to both sides
and then using the Taylor theorem we have
\begin{eqnarray}
(x_{0}-x_{2})^{k}x_{0}^{l}Y_{W}(u,x_{0})Y_{W}(v,x_{2})w
=e^{-x_{2}{\partial\over \partial x_{0}}}p(x_{0},x_{2})
=p(x_{0}-x_{2},x_{2}).
\end{eqnarray}
Hence
\begin{eqnarray}
(x_{0}-x_{2})^{k}Y_{W}(u,x_{0})Y_{W}(v,x_{2})w=x_{0}^{-l}p(x_{0}-x_{2},x_{2})
\in W((x_{0},x_{2})).
\end{eqnarray}
Since $k$ is independent of $w$, (\ref{ecompat2vo}) follows.
$\;\;\;\;\Box$

\br{rll}
{\em Lemma \ref{lcompat2vo}
in a slightly different form has been obtained
in [LL] (Proposition 3.3.12).}
\er

To generalize  Lemma \ref{lcompat2vo} for the products of more than two vertex
operators we shall need to assume that $V$ is an axiomatic 
$G_{1}$-vertex algebra and assume that the stronger weak associativity
holds on $W$.

\bp{pcompatibilitycon}
Let $V$ be an axiomatic $G_{1}$-vertex algebra and 
let $(W,Y_{W})$ be a $V$-module (with the stronger weak associativity).
Then for any $v^{(1)},\dots, v^{(r)}\in V$,
there exists a nonnegative integer $k$ such that
\begin{eqnarray}
\left(\prod_{1\le i<j\le r}(x_{i}-x_{j})^{k}\right)
Y_{W}(v^{(1)},x_{1})\cdots Y_{W}(v^{(r)},x_{r})
\in \Hom (W,W((x_{1},\dots,x_{r}))).
\end{eqnarray}
\ep

{\bf Proof.} First, we prove the special case with $W=V$
by induction on $r$.
For $r=2$, it has been proved by Lemma \ref{lcompat2vo}. 
Assume that the assertion holds for 
a certain $r\ge 2$. Let
$v^{(1)},\dots, v^{(r+1)}\in V$ and $w\in W$.
{}From the (stronger) weak associativity, 
there exists $l\in \N$ such that 
\begin{eqnarray}
(x_{0i}+x_{r+1})^{l}Y(v^{(i)},x_{0i}+x_{r+1})Y(v,x_{r+1})w
=(x_{0i}+x_{r+1})^{l}Y(Y(v^{(i)},x_{0i})v,x_{r+1})w
\end{eqnarray}
for {\em all} $v\in V$ and for $i=1,\dots,r$.
Then
\begin{eqnarray}\label{egeneralcomalg}
& &\left(\prod_{i=1}^{r}(x_{0i}+x_{r+1})^{l}\right)
Y(v^{(1)},x_{01}+x_{r+1})\cdots Y(v^{(r)},x_{0r}+x_{r+1})
Y(v^{(r+1)},x_{r+1})w\nonumber\\
&=&\left(\prod_{i=1}^{r}(x_{0i}+x_{r+1})^{l}\right)
Y(v^{(1)},x_{01}+x_{r+1})\cdots Y(v^{(r-1)},x_{0r-1}+x_{r+1})\nonumber\\
& &\cdot Y(Y(v^{(r)},x_{0r})v^{(r+1)},x_{r+1})w\nonumber\\
&=&\left(\prod_{i=1}^{r}(x_{0i}+x_{r+1})^{l}\right)
Y(Y(v^{(1)},x_{01})\cdots Y(v^{(r)},x_{0r})v^{(r+1)},x_{r+1})w.
\end{eqnarray}
Notice that for the second equality we are using the stronger
version of the weak associativity.
{}From the inductive hypothesis there exists a nonnegative integer $k'$
such that
\begin{eqnarray}\label{egeneralcomalgind0}
\left(\prod_{1\le i<j\le r}(x_{0i}-x_{0j})^{k'}\right)
Y(v^{(1)},x_{01})\cdots Y(v^{(r)},x_{0r})v^{(r+1)}\in V((x_{01},\dots,x_{0r})),
\end{eqnarray}
so that there exists a nonnegative integer $k''$
(only depending on $v^{(i)}$'s, not $w$) such that
\begin{eqnarray}\label{egeneralcomalgind}
\left(\prod_{1\le i<j\le r}
(x_{0i}-x_{0j})^{k'}\right)x_{01}^{k''}\cdots x_{0r}^{k''}
Y(v^{(1)},x_{01})\cdots Y(v^{(r)},x_{0r})v^{(r+1)}\in V[[x_{01},\dots,x_{0r}]].
\end{eqnarray}
Combining (\ref{egeneralcomalg}) with (\ref{egeneralcomalgind}) we get
\begin{eqnarray}
& &\left(\prod_{i=1}^{r}(x_{0i}+x_{r+1})^{l}\prod_{1\le i<j\le r}
(x_{0i}-x_{0j})^{k'}\right)x_{01}^{k''}\cdots x_{0r}^{k''}\cdot \nonumber\\
& &\cdot Y(v^{(1)},x_{01}+x_{r+1})\cdots Y(v^{(r)},x_{0r}+x_{r+1})
Y(v^{(r+1)},x_{r+1})w\nonumber\\
& &\in W[[x_{01},\dots,x_{0r}]]((x_{r+1})).
\end{eqnarray}
Therefore (by substituting $x_{0i}=x_{i}-x_{r+1}$)
\begin{eqnarray}
& &\left(\prod_{i=1}^{r}x_{i}^{l}\prod_{1\le i<j\le r}
(x_{i}-x_{j})^{k'}\right)(x_{1}-x_{r+1})^{k''}\cdots
(x_{r}-x_{r+1})^{k''}\cdot \nonumber\\
& &\cdot Y(v^{(1)},x_{1})\cdots Y(v^{(r)},x_{r})
Y(v^{(r+1)},x_{r+1})w\nonumber\\
& &\in W[[x_{1},\dots,x_{r}]]((x_{r+1})).
\end{eqnarray}
That is,
\begin{eqnarray}
& &\left(\prod_{1\le i<j\le r}
(x_{i}-x_{j})^{k'}\right)(x_{1}-x_{r+1})^{k''}\cdots
(x_{r}-x_{r+1})^{k''}\cdot \nonumber\\
& &\cdot Y(v^{(1)},x_{1})\cdots Y(v^{(r)},x_{r})
Y(v^{(r+1)},x_{r+1})w\nonumber\\
& &\in W((x_{1},\dots,x_{r},x_{r+1})).
\end{eqnarray}
Since $k'$ and $k''$ are independent of $w$, this 
finishes the induction, proving the special case.

For the general case, similar to the special case,
(from the (stronger) weak associativity), there exists $l\in \N$ 
such that
\begin{eqnarray}\label{generalcommod}
& &\left(\prod_{i=1}^{r}(x_{0i}+x_{r+1})^{l}\right)
Y_{W}(v^{(1)},x_{01}+x_{r+1})\cdots Y_{W}(v^{(r)},x_{0r}+x_{r+1})
Y_{W}(v^{(r+1)},x_{r+1})w\nonumber\\
&=&\left(\prod_{i=1}^{r}(x_{0i}+x_{r+1})^{l}\right)
Y_{W}(Y(v^{(1)},x_{01})\cdots Y(v^{(r)},x_{0r})v^{(r+1)},x_{r+1})w.
\end{eqnarray}
The rest directly follows from 
the proof and the result of the special case.
$\;\;\;\;\Box$

\br{rrationality}
{\em Let $V$ and $W$ be given as in Proposition \ref{pcompatibilitycon}.
Let $w^{*}\in W^{*},\; v^{(1)},\dots,v^{(r)}\in V,\; w\in W$. 
In view of Proposition \ref{pcompatibilitycon},
there exist nonnegative integers $k$ and $l$ such that
\begin{eqnarray}
\<w^{*},Y_{W}(v^{(1)},x_{1})\cdots Y_{W}(v^{(r)},x_{r})w\>
=\left(\prod_{1\le i<j\le r}(x_{i}-x_{j})^{-k}\prod_{i=1}^{r}x_{i}^{-l}
\right)p(x_{1},\dots,x_{r})
\end{eqnarray}
for some $p(x_{1},\dots,x_{r})\in \C[[x_{1},\dots,x_{r}]]$, where
we are using the binomial expansion convention.}
\er

\br{rsimplicity}
{\em For an axiomatic $G_{1}$-vertex algebra $V$,
a left ideal of $V$ amounts to a $V$-submodule of $V$.
Unlike in the case of ordinary vertex algebras, 
the notions of left and right ideals are in general different. 
Thus, the simplicity of $V$ as a $V$-module
does not amount to the simplicity of $V$ as an axiomatic 
$G_{1}$-vertex algebra.}
\er

\section{Axiomatic $G_{1}$-vertex algebras generated 
by compatible $G_{1}$-vertex operators}

In this section we study (weak) $G_{1}$-vertex operators on an arbitrary
vector space and we show how a suitable set of (weak) $G_{1}$-vertex operators 
gives rise to an axiomatic $G_{1}$-vertex algebra. We recover
the corresponding result of [Li1].

Let $W$ be a vector space fixed throughout this section. 

\bd{dg1vertexoperator}
{\em A {\em weak $G_{1}$-vertex operator} on $W$ 
is a formal series 
\begin{eqnarray}
a(x)=\sum_{n\in \Z}a_{n}x^{-n-1}\in (\End W)[[x,x^{-1}]]
\end{eqnarray}
such that for every $w\in W$, $a_{n}w=0$ for $n$ sufficiently
large. Namely, a weak $G_{1}$-vertex operator on $W$ 
is an element of $\Hom (W,W((x)))$.}
\ed

All weak $G_{1}$-vertex operators on $W$
constitute the space $\Hom (W,W((x)))$.
We alternatively denote this space by ${\cal{E}}_{G_{1}}(W)$.

\br{rli2weakvo}
{\em In [Li1] and [LL], an element of $\Hom (W,W((x)))$ was simply called
a weak vertex operator. We here chose to use the term 
``weak $G_{1}$-vertex operator''
because of the expected connection 
with Borcherds' notion of $G$-vertex algebra.} 
\er

Set
\begin{eqnarray}
D={d\over dx}.
\end{eqnarray}
Then $D$ is a natural endomorphism of ${\cal{E}}_{G_{1}}(W)$.

Motivated by Proposition  \ref{pcompatibilitycon} and 
by [B2] we define the following notion of compatibility.

\bd{dcompatibility}
{\em An (ordered) sequence $(\psi^{(1)},\dots,\psi^{(r)})$ 
in ${\cal{E}}_{G_{1}}(W)$ is said to be {\em compatible} 
if there exists a nonnegative integer $k$ such that}
\begin{eqnarray}
\left(\prod_{1\le i<j\le r}(x_{i}-x_{j})^{k}\right)
\psi^{(1)}(x_{1})\cdots \psi^{(r)}(x_{r})
\in \Hom (W,W((x_{1},\dots,x_{r}))).
\end{eqnarray}
\ed

A set or a space $S$ of weak $G_{1}$-vertex operators on $W$ is said 
to be compatible if any finite sequence in $S$ is compatible.
A weak $G_{1}$-vertex operator $a(x)$ on $W$ is called a {\em
$G_{1}$-vertex operator} if $\{ a(x)\}$ is compatible.
Then weak $G_{1}$-vertex operators in a compatible set are
$G_{1}$-vertex operators.
It is important to note that compatibility in general 
depends on the order. Clearly, $(\End W)((x))$ is
a compatible subspace of ${\cal{E}}_{G_{1}}(W)$. 

\br{rspan}
{\em It is easy to see that the linear span of any compatible set
of weak $G_{1}$-vertex operators on $W$ is compatible.}
\er

\bex{exacompatibleset}
{\em Let $V$ be an axiomatic $G_{1}$-vertex algebra and 
let $(W,Y_{W})$ be a $V$-module. It follows immediately from 
Proposition \ref{pcompatibilitycon} that
the image of $V$ under $Y_{W}$ is a compatible
space of weak $G_{1}$-vertex operators on $W$.}
\eex

It is in general not a good idea
to use the definition directly to check the compatibility of 
a set of weak
$G_{1}$-vertex operators. In the following, we prove that
certain pairwise relations imply compatibility.

\bl{lgeneralizedcomm}
Let $a(x), b(x)\in {\cal{E}}_{G_{1}}(W)$. Assume that there exists
a nonnegative integer $k$ such that
\begin{eqnarray}\label{ehalflocality}
(x_{1}-x_{2})^{k}a(x_{1})b(x_{2})=(x_{1}-x_{2})^{k}\sum_{i=1}^{r}
\psi^{(i)}(x_{2})\phi^{(i)}(x_{1})
\end{eqnarray}
for some $\psi^{(i)}(x),\phi^{(i)}(x)\in {\cal{E}}_{G_{1}}(W)$.
Then the ordered sequence $(a(x), b(x))$ is compatible.
\el

{\pf} Let $w\in W$. Since $b(x_{2})\in \Hom (W,W((x_{2})))$, 
$(x_{1}-x_{2})^{k}a(x_{1})b(x_{2})w$ involves only finitely many negative 
powers of $x_{2}$. On the other hand, the expression on the right-hand side of
 (\ref{ehalflocality}), after applied to $w$, involves only finitely many 
negative powers of $x_{1}$. Consequently, 
$$(x_{1}-x_{2})^{k}a(x_{1})b(x_{2})w\in W((x_{1},x_{2})).$$
Thus $(x_{1}-x_{2})^{k}a(x_{1})b(x_{2})\in \Hom (W,W((x_{1},x_{2})))$.
That is, $(a(x),b(x))$ is compatible.
$\;\;\;\;\Box$

Furthermore, we have:

\bp{ppracticalcase}
Let $S$ be a set of weak $G_{1}$-vertex operators on $W$ such that
for any $a(x),b(x)\in S$, there exists a nonnegative integer
$k$ such that 
\begin{eqnarray}\label{e4.3}
(x_{1}-x_{2})^{k}a(x_{1})b(x_{2})
=\sum_{i=1}^{r}\alpha_{i} (x_{1}-x_{2})^{k}b^{(i)}(x_{2})a^{(i)}(x_{1})
\end{eqnarray}
for some $\alpha_{i}\in \C,\; a^{(i)}(x), b^{(i)}(x)\in S$.
Then $S$ is compatible.
\ep

{\pf} We must prove that any sequence in $S$ of finite length
is compatible. We shall use induction on the length $n$ of sequences.
If $n=2$, it has been proved by Lemma \ref{lgeneralizedcomm}.
Assume that any sequence  in $S$ of length $n$ is compatible.
Let $\psi^{(1)},\dots,\psi^{(n+1)}\in S$. From the inductive hypothesis,
there exists a nonnegative integer $k_{1}$ such that
\begin{eqnarray}\label{eproduct1}
\left(\prod_{2\le i<j\le n+1}(x_{i}-x_{j})^{k_{1}}\right)
\psi^{(2)}(x_{2})\cdots \psi^{(n+1)}(x_{n+1})
\in \Hom (W,W((x_{2},\dots,x_{n+2}))).
\end{eqnarray}
{}From (\ref{e4.3}) there exists a nonnegative integer $k_{2}$ such that
\begin{eqnarray}\label{eproduct2}
(x_{1}-x_{2})^{k_{2}}\psi^{(1)}(x_{1})\psi^{(2)}(x_{2})
=\sum_{i=1}^{r}\alpha_{i}(x_{1}-x_{2})^{k_{2}}b^{(i)}(x_{2})a^{(i)}(x_{1})
\end{eqnarray}
for some $\alpha_{i}\in \C,\; a^{(i)}(x), b^{(i)}(x)\in S$.
{}From the inductive hypothesis again, there exists 
a nonnegative integer $k_{3}$ such that
\begin{eqnarray}\label{eproduct3}
& &\left(\prod_{i=3}^{n+1}(x_{1}-x_{i})^{k_{3}}
\prod_{3\le i<j\le n+1}(x_{i}-x_{j})^{k_{3}}\right)
a^{(s)}(x_{1})\psi^{(3)}(x_{3})\cdots \psi^{(n+1)}(x_{n+1})
\nonumber\\
&\in& \Hom (W, W((x_{1},x_{3},x_{4},\dots,x_{n+1})))
\end{eqnarray}
for $s=1,\dots,r$. Because of (\ref{eproduct2}) we have
\begin{eqnarray}\label{eproduct4}
& &(x_{1}-x_{2})^{k_{2}}\left(\prod_{i=3}^{n+1}(x_{1}-x_{i})^{k_{3}}
\prod_{3\le i<j\le n+1}(x_{i}-x_{j})^{k_{3}}\right)
\psi^{(1)}(x_{1})\cdots \psi^{(n+1)}(x_{n+1})\nonumber\\
&=&
(x_{1}-x_{2})^{k_{2}}\left(\prod_{i=3}^{n+1}(x_{1}-x_{i})^{k_{3}}
\prod_{3\le i<j\le n+1}(x_{i}-x_{j})^{k_{3}}\right)\nonumber\\
& &\cdot \sum_{s=1}^{r}\alpha_{i}
b^{(s)}(x_{2})a^{(s)}(x_{1})\psi^{(3)}(x_{3})\cdots \psi^{(n+1)}(x_{n+1}).
\end{eqnarray}
{}From (\ref{eproduct3}), the right-hand side of (\ref{eproduct4}) lies in
$$\Hom (W,W((x_{2}))((x_{1},x_{3},x_{4},\dots,x_{n+1}))),$$
and so does the left-hand side of (\ref{eproduct4}).
Combining this with (\ref{eproduct1}) we see that
\begin{eqnarray}
& &\left(\prod_{1\le i<j\le n+1}(x_{i}-x_{j})^{k_{1}+k_{2}+k_{3}}\right)
\psi^{(1)}(x_{1})\psi^{(2)}(x_{2})\cdots \psi^{(n+1)}(x_{n+1})\nonumber\\
&\in &\Hom (W,W((x_{1},x_{2},\dots,x_{n+1}))).
\end{eqnarray}
This proves that the sequence $(\psi^{(1)},\dots,\psi^{(n+1)})$ 
is compatible, completing the induction.
$\;\;\;\;\Box$

\br{rlocality}
{\em Recall from [Li1] that weak ($G_{1}$-)vertex operators 
$a(x)$ and $b(x)$ are said to be {\em mutually local} 
if there exists a nonnegative integer $k$ such that
\begin{eqnarray}
(x_{1}-x_{2})^{k}a(x_{1})b(x_{2})
=(x_{1}-x_{2})^{k}b(x_{2})a(x_{1}).
\end{eqnarray}
A set $S$ of weak vertex operators on $W$ is said to be local if 
any two (maybe the same) weak vertex operators in $S$ are mutually local.}
\er

As an immediate consequence of Proposition \ref{ppracticalcase}
we have:

\bc{clocalimplycompatibility}
Any local set of weak $G_{1}$-vertex operators on $W$ is compatible.
$\;\;\;\;\Box$
\ec

\bl{lpreparefordefinition}
Let $a(x), b(x)\in {\cal{E}}_{G_{1}}(W)$ be such that 
the sequence $(a(x),b(x))$ is compatible.
Define
\begin{eqnarray}
T(a(x)b(y))=(-y+x)^{-k}\left((x-y)^{k}a(x)b(y)\right),
\end{eqnarray}
where $k$ is any nonnegative integer such that
\begin{eqnarray}\label{exykab}
(x-y)^{k}a(x)b(y)\in \Hom(W,W((x,y))).
\end{eqnarray}
Then $T(a(x)b(y))$ does not depend the choice of $k$
and it lies in $\Hom (W,W((y))((x)))$.
Furthermore, 
\begin{eqnarray}\label{etab=ab}
(x-y)^{k}T(a(x)b(y))=(x-y)^{k}a(x)b(y)
\end{eqnarray}
for any nonnegative integer $k$ such that (\ref{exykab}) holds.
\el

{\bf Proof.} Clearly, for any nonnegative integer $k$ such that 
(\ref{exykab}) holds, we have
\begin{eqnarray}
(-y+x)^{-k}\left((x-y)^{k}a(x)b(y)\right)\in \Hom (W,W((y))((x))).
\end{eqnarray}
We easily see that (\ref{etab=ab}) holds if $T(a(x)b(y))$ is well defined.
It remains to prove that $T(a(x)b(y))$ is independent of $k$.

Let $k_{1}$ and $k_{2}$ be any two nonnegative integers such that
\begin{eqnarray}
(x-y)^{k_{i}}a(x)b(y)\in \Hom(W,W((x,y)))\;\;\;\mbox{ for }i=1,2.
\end{eqnarray}
Assume that $k_{1}\ge k_{2}$.
(For the case that $k_{2}\ge k_{1}$, one simply 
exchanges $k_{1}$ with $k_{2}$ in the following argument.)
Then
\begin{eqnarray}
(-y+x)^{-k_{1}}\left((x-y)^{k_{1}}a(x)b(y)\right)
&=&(-y+x)^{-k_{1}}\left((x-y)^{k_{1}-k_{2}}\left((x-y)^{k_{2}}a(x)b(y)\right)
\right)\nonumber\\
&=&(-y+x)^{-k_{1}}(x-y)^{k_{1}-k_{2}}\left((x-y)^{k_{2}}a(x)b(y)\right)
\nonumber\\
&=&(-y+x)^{-k_{2}}\left((x-y)^{k_{2}}a(x)b(y)\right).
\end{eqnarray}
This proves the assertion. $\;\;\;\;\Box$

\bd{danbdefinition}
{\em Let $a(x), b(x)\in {\cal{E}}_{G_{1}}(W)$ be such that 
the sequence $(a(x),b(x))$ is compatible. For $n\in \Z$, we define
$a(x)_{n}b(x)\in (\End W)[[x,x^{-1}]]$ by
\begin{eqnarray}
a(x)_{n}b(x)=\Res_{x_{1}}\left( (x_{1}-x)^{n}a(x_{1})b(x)
-(-x+x_{1})^{n}T(a(x_{1})b(x))\right).
\end{eqnarray}}
\ed

Just as in [Li1] and [LL] with Lemma \ref{lpreparefordefinition}
 we immediately have:

\bp{pnthproduct}
Let $a(x), b(x)\in {\cal{E}}_{G_{1}}(W)$ be such that 
the sequence $(a(x),b(x))$ is compatible. We have
\begin{eqnarray}
a(x)_{n}b(x)\in {\cal{E}}_{G_{1}}(W)\;\;\;\mbox{ for }n\in \Z.
\end{eqnarray}
Furthermore,
\begin{eqnarray}\label{etruncationresult}
a(x)_{n}b(x)=0\;\;\;\mbox{ for }n\ge k,
\end{eqnarray}
where $k$ is a nonnegative integer such that 
$(x-y)^{k}a(x)b(y)\in \Hom (W,W((x,y)))$.$\;\;\;\;\;\Box$
\ep

\br{radvantage}
{\em Let $\psi, \phi, \psi^{(i)}, \phi^{(i)}\in {\cal{E}}_{G_{1}}(W)$ 
for $i=1,\dots,r$ be such that
\begin{eqnarray}
(x_{1}-x_{2})^{k}\psi(x_{1})\phi(x_{2})
=(x_{1}-x_{2})^{k}\sum_{i=1}^{r}\alpha_{i}\phi^{(i)}(x_{2})\psi^{(i)}(x_{1})
\end{eqnarray}
for some nonnegative integer $k$ and some $\alpha_{i}\in \C$. 
(In view of Lemma \ref{lgeneralizedcomm}, the sequence $(\psi,\phi)$ 
is compatible.) Then
\begin{eqnarray}
T(\psi(x_{1})\phi(x_{2}))
&=&(-x_{2}+x_{1})^{-k}\left((x_{1}-x_{2})^{k}\psi(x_{1})\phi(x_{2})\right)
\nonumber\\
&=&(-x_{2}+x_{1})^{-k}\left((x_{1}-x_{2})^{k}\sum_{i=1}^{r}
\alpha_{i}\phi^{(i)}(x_{2})\psi^{(i)}(x_{1})\right)\nonumber\\
&=&\sum_{i=1}^{r}
\alpha_{i}\phi^{(i)}(x_{2})\psi^{(i)}(x_{1}).
\end{eqnarray}
Therefore,
\begin{eqnarray}
\psi(x)_{n}\phi(x)&=&\Res_{x_{1}}\left( (x_{1}-x)^{n}\psi(x_{1})\phi(x)
-\sum_{i=1}^{r}\alpha_{i} (-x+x_{1})^{n}\phi^{(i)}(x)\psi^{(i)}(x_{1})\right).
\end{eqnarray}
In particular, if
$$(x_{1}-x_{2})^{k}\psi(x_{1})\phi(x_{2})
=\alpha (x_{1}-x_{2})^{k}\phi(x_{2})\psi(x_{1})$$
for some $\alpha\in \C$, we have 
\begin{eqnarray}
\psi(x)_{n}\phi(x)&=&\Res_{x_{1}}\left( (x_{1}-x)^{n}\psi(x_{1})\phi(x)
-\alpha (-x+x_{1})^{n}\phi(x)\psi(x_{1})\right).
\end{eqnarray}}
\er

\br{radvantageanddisadvantage}
{\em In view of Remark \ref{radvantage},
if $a(x), b(x)$ are mutually local weak vertex operators on $W$,
then the current definition for $a(x)_{n}b(x)$ 
coincides with the one given in [Li1] and [LL],
where for {\em any} $\alpha(x), \beta(x)\in {\cal{E}}_{G_{1}}(W)$
it was defined that
\begin{eqnarray*}
\alpha(x)_{n}\beta(x)
=\Res_{x_{1}}\left((x_{1}-x)^{n}\alpha(x_{1})\beta(x)
-(-x+x_{1})^{n}\beta(x)\alpha(x_{1})\right).
\end{eqnarray*}
If $a(x),b(x)$ are weak vertex operators with the relation
$$(x_{1}-x_{2})^{k}a(x_{1})b(x_{2})=-(x_{1}-x_{2})^{k}b(x_{2})a(x_{1}),$$
then the current definition for $a(x)_{n}b(x)$ 
is different from the one given in [Li1] and [LL].
The essential difference between the definitions is that 
the current definition only uses the product $a(x_{1})b(x_{2})$,
not the product $b(x_{2})a(x_{1})$ while the definition 
given in [Li1] and [LL] uses both of the products. That is, 
one definition takes the associative
algebra point of view and the other takes the Lie algebra 
point of view.}
\er

Writing $a(x)_{n}b(x)$ for $n\in \Z$ in terms of generating function as
\begin{eqnarray}
Y_{\cal{E}}(a(x),x_{0})b(x)
=\sum_{n\in \Z}a(x)_{n}b(x) x_{0}^{-n-1},
\end{eqnarray}
we have
\begin{eqnarray}
Y_{\cal{E}}(a(x),x_{0})b(x)\in {\cal{E}}_{G_{1}}(W)((x_{0})).
\end{eqnarray}
Then
\begin{eqnarray}\label{ecanonicalmodulejacobi}
& &Y_{\cal{E}}(a(x_{2}),x_{0})b(x_{2})
=\sum_{n\in \Z}(a(x_{2})_{n}b(x_{2}))x_{0}^{-n-1}\nonumber\\
&=&\Res_{x_{1}}\left
( x_{0}^{-1}\delta\left(\frac{x_{1}-x_{2}}{x_{0}}\right)
a(x_{1})b(x_{2})
-x_{0}^{-1}\delta\left(\frac{x_{2}-x_{1}}{-x_{0}}\right)
T(a(x_{1})b(x_{2}))\right).
\end{eqnarray}

\br{rdomainYE}
{\em Since $a(x)_{n}b(x)$ for $n\in \Z$ are defined 
under the condition that the (ordered) sequence $(a(x),b(x))$ 
is compatible, $Y_{\cal{E}}(a(x),x_{0})b(x)$ is a well defined 
element of ${\cal{E}}_{G_{1}}(W)((x_{0}))$ 
under the same condition. 
For any compatible space $V$ of $G_{1}$-vertex operators on $W$,
$Y_{\cal{E}}$ is a natural linear map from $V$ 
to $\Hom (V,{\cal{E}}_{G_{1}}(W)((x_{0})))$.}
%Is there a canonical linear extension of $Y_{\cal{E}}$?}
\er

\bp{pdefinitionassociativity}
Let $a(x), b(x)\in {\cal{E}}_{G_{1}}(W)$ be such that 
the sequence $(a(x),b(x))$ is compatible.
Then for any $w\in W$, there exists a nonnegative integer
$l$ such that
\begin{eqnarray}\label{e02labw}
(x_{0}+x_{2})^{l}a(x_{0}+x_{2})b(x_{2})w\in W((x_{0},x_{2})).
\end{eqnarray}
Furthermore, if $l$ is a nonnegative integer such that (\ref{e02labw})
holds, then
\begin{eqnarray}\label{edefinitionassociativity}
(x_{2}+x_{0})^{l}(Y_{\cal{E}}(a(x_{2}),x_{0})b(x_{2}))w=
(x_{0}+x_{2})^{l}a(x_{0}+x_{2})b(x_{2})w.
\end{eqnarray}
\ep

{\bf Proof.} Let $k$ be a nonnegative integer such that
$$(x_{1}-x_{2})^{k}a(x_{1})b(x_{2})\in \Hom (W,W((x_{1},x_{2}))).$$
In view of Lemma \ref{lpreparefordefinition}, we have
$$(x_{1}-x_{2})^{k}T(a(x_{1})b(x_{2}))=(x_{1}-x_{2})^{k}a(x_{1})b(x_{2}).$$
For any $w\in W$, since 
$(x_{1}-x_{2})^{k}a(x_{1})b(x_{2})w\in W((x_{1},x_{2}))$, there exists 
a nonnegative integer $l$ such that 
$$x_{1}^{l}(x_{1}-x_{2})^{k}a(x_{1})b(x_{2})w
\in W[[x_{1},x_{2}]][x_{2}^{-1}].$$
Then
\begin{eqnarray}
(x_{0}+x_{2})^{l}x_{0}^{k}a(x_{0}+x_{2})b(x_{2})w
\in W[[x_{0},x_{2}]][x_{2}^{-1}],
\end{eqnarray}
which implies (\ref{e02labw}).

Now let $l$ be a nonnegative integer such that (\ref{e02labw}) holds.
Let $k'$ be a nonnegative integer such that
$$x_{0}^{k'}(x_{0}+x_{2})^{l}a(x_{0}+x_{2})b(x_{2})w
\in W[[x_{0},x_{2}]][x_{2}^{-1}].$$
Then
\begin{eqnarray}\label{exlxykabw}
x_{1}^{l}(x_{1}-x_{2})^{k'}a(x_{1})b(x_{2})w\in W[[x_{1},x_{2}]][x_{2}^{-1}].
\end{eqnarray}
Therefore,
\begin{eqnarray}
x_{1}^{l}(x_{1}-x_{2})^{k+k'}T(a(x_{1})b(x_{2}))w
=x_{1}^{l}(x_{1}-x_{2})^{k+k'}a(x_{1})b(x_{2})w
\in W[[x_{1},x_{2}]][x_{2}^{-1}].
\end{eqnarray}
Multiplying by $(-x_{2}+x_{1})^{-k-k'}$, which lies 
in $\C[x_{2},x_{2}^{-1}][[x_{1}]]$, we get
\begin{eqnarray}
x_{1}^{l}T(a(x_{1})b(x_{2}))w\in W((x_{2}))[[x_{1}]].
\end{eqnarray}
Just as in the ordinary vertex algebra theory (cf. [DL], [Li1] or [LL]), 
multiplying both sides of (\ref{ecanonicalmodulejacobi})
by $(x_{2}+x_{0})^{l}$ we obtain (\ref{edefinitionassociativity}).
$\;\;\;\;\Box$

\br{ralternativedef}
{\em In view of (\ref{etruncationresult}), by
multiplying both sides of (\ref{edefinitionassociativity})
by $(x_{2}+x_{0})^{-1}$ we get
\begin{eqnarray}\label{ealternativedef}
(Y_{\cal{E}}(a(x_{2}),x_{0})b(x_{2}))w=(x_{2}+x_{0})^{-l}
\left[(x_{0}+x_{2})^{l}a(x_{0}+x_{2})b(x_{2})w\right].
\end{eqnarray}
On the other hand, notice that (\ref{exlxykabw}) implies that 
\begin{eqnarray}\label{exylabw}
(x_{0}+x_{2})^{l}a(x_{0}+x_{2})b(x_{2})w\in W((x_{0},x_{2})),
\end{eqnarray}
so that the expression on the right-hand side of 
(\ref{ealternativedef}) is well defined.
Then one can use (\ref{exlxykabw}) as an alternative definition for 
$Y_{\cal{E}}(a(x_{2}),x_{0})b(x_{2})$.}
\er

\br{rjacobidef}
{\em Combining Lemma \ref{lpreparefordefinition} and
Proposition \ref{pdefinitionassociativity} with
Lemma \ref{lformaljacobiidentity} we get
\begin{eqnarray}\label{emoduleactionjacobi}
& &x_{2}^{-1}\delta\left(\frac{x_{1}-x_{0}}{x_{2}}\right)
Y_{\cal{E}}(a(x_{2}),x_{0})b(x_{2})\nonumber\\
&=&x_{0}^{-1}\delta\left(\frac{x_{1}-x_{2}}{x_{0}}\right)a(x_{1})b(x_{2})
-x_{0}^{-1}\delta\left(\frac{x_{2}-x_{1}}{-x_{0}}\right)
T(a(x_{1})b(x_{2})).
\end{eqnarray}}
\er

Let $1_{W}$ denote the identity operator on $W$ and let
$a(x)$ be any weak $G_{1}$-vertex operator on $W$. 
Since
\begin{eqnarray}
1_{W}(x_{1})a(x_{2})&=&a(x_{2})\in \Hom (W,W((x_{2})))\subset
\Hom (W,W((x_{1},x_{2}))),\\
a(x_{1})1_{W}(x_{2})&=&a(x_{1})\in \Hom (W,W((x_{1})))\subset
\Hom (W,W((x_{1},x_{2}))),
\end{eqnarray}
the sequences $(1_{W},a(x))$ and $(a(x),1_{W})$ are compatible, 
 and $T(a(x_{1})1_{W})=a(x_{1})$.
By (\ref{ecanonicalmodulejacobi}) we have
\begin{eqnarray}
& &Y_{\cal{E}}(1_{W}(x),x_{0})a(x)\nonumber\\
&=&\Res_{x_{1}}\left(x_{0}^{-1}\delta\left(\frac{x_{1}-x}{x_{0}}\right)
1_{W}(x_{1})a(x)-x_{0}^{-1}\delta\left(\frac{x-x_{1}}{-x_{0}}\right)
a(x)\right)\nonumber\\
&=&\Res_{x_{1}}x^{-1}\delta\left(\frac{x_{1}-x_{0}}{x}\right)
a(x)\nonumber\\
&=&a(x).
\end{eqnarray}
Similarly we have
\begin{eqnarray}
Y_{\cal{E}}(a(x),x_{0})1_{W}(x)
&=&\Res_{x_{1}}x^{-1}\delta\left(\frac{x_{1}-x_{0}}{x}\right)a(x_{1})
\nonumber\\
&=&a(x+x_{0})\nonumber\\
&=&e^{x_{0}{d\over d x}}a(x)\nonumber\\
&=&e^{x_{0}D}a(x).
\end{eqnarray}
Thus we have proved:

\bl{lvacuum-creationproperties}
For any $a(x)\in {\cal{E}}_{G_{1}}(W)$,
\begin{eqnarray}
& &Y_{\cal{E}}(1_{W},x_{0})a(x)=a(x),\\
& &Y_{\cal{E}}(a(x),x_{0})1_{W}=e^{x_{0}D}a(x).\;\;\;\;\;\Box
\end{eqnarray}
\el

A compatible space $U$ of weak $G_{1}$-vertex operators on $W$
is said to be {\em closed} if 
\begin{eqnarray}
a(x)_{n}b(x)\in U\;\;\;\mbox{ for }a(x),b(x)\in
U,\; n\in \Z.
\end{eqnarray}
Then for a closed compatible space $U$, we have a linear map
$Y_{\cal{E}}$ from $U$ to $\Hom (U, U((x_{0})))$.

\br{rclosedcompatibleset}
{\em Let $V$ be an axiomatic $G_{1}$-vertex algebra and let $(W,Y_{W})$
be a $V$-module. Then the image of $Y_{W}$ is a closed compatible subspace
of ${\cal{E}}_{G_{1}}(W)$. Furthermore, for $u,v\in V,\; n\in \Z$,
\begin{eqnarray}\label{eywunywv}
Y_{W}(u_{n}v,x)=Y_{W}(u,x)_{n}Y_{W}(v,x).
\end{eqnarray}
Indeed, for any $u,v\in V,\; w\in W$, in view of Definition
\ref{dmodule} and
Proposition \ref{pdefinitionassociativity}, there exists a nonnegative
integer $l$ such that
\begin{eqnarray}
(x_{0}+x_{2})^{l}Y_{W}(Y(u,x_{0})v,x_{2})w
&=&(x_{0}+x_{2})^{l}Y_{W}(u,x_{0}+x_{2})Y_{W}(v,x_{2})w\nonumber\\
(x_{0}+x_{2})^{l}(Y_{\cal{E}}(Y_{W}(u,x_{2}),x_{0})Y_{W}(v,x_{2}))w
&=&(x_{0}+x_{2})^{l}Y_{W}(u,x_{0}+x_{2})Y_{W}(v,x_{2})w.
\end{eqnarray}
Then
\begin{eqnarray}
(x_{0}+x_{2})^{l}Y_{W}(Y(u,x_{0})v,x_{2})w
=(x_{0}+x_{2})^{l}(Y_{\cal{E}}(Y_{W}(u,x_{2}),x_{0})Y_{W}(v,x_{2}))w.
\end{eqnarray}
Noticing that both sides involve only finitely many negative powers of
$x_{0}$, by multiplying by $(x_{2}+x_{0})^{-l}$ we get
\begin{eqnarray}
Y_{W}(Y(u,x_{0})v,x_{2})w=(Y_{\cal{E}}(Y_{W}(u,x_{2}),x_{0})Y_{W}(v,x_{2}))w,
\end{eqnarray}
which is (\ref{eywunywv}) in terms of generating functions.}
\er

In view of Remark \ref{rclosedcompatibleset}, if $(W,Y_{W})$ 
is a faithful $V$-module, e.g., $W=V$ (the faithfulness follows from
the creation property), $V$ can be naturally
identified with a closed compatible subspace of
${\cal{E}}_{G_{1}}(W)$, containing $1_{W}$.
Next, we shall show that for an abstract vector space $W$,
any closed compatible subspace of ${\cal{E}}_{G_{1}}(W)$ that contains $1_{W}$
is a (weak) axiomatic $G_{1}$-vertex algebra with 
$W$ as a natural faithful module.

We here introduce a notation for convenience.
Let $U$ be a vector space and let
$$a(x)=\sum_{n\in \Z}a_{n}x^{-n-1}\in U[[x,x^{-1}]]$$ 
be any formal series, e.g., a weak $G_{1}$-vertex
operator on $W$.  For $m\in \Z$, we set
\begin{eqnarray}
a(x)_{\ge m}=\sum_{n\ge m}a_{n}x^{-n-1}.
\end{eqnarray}
Then for any polynomial $p(x)$ we have
\begin{eqnarray}\label{esimplefactresidule}
\Res_{x}x^{m}p(x)a(x)=\Res_{x}x^{m}p(x)a(x)_{\ge m}.
\end{eqnarray}
First, we have the following result:

\bl{lbeforegenerating}
Let $\psi_{1},\dots, \psi_{r}, a,b, \phi_{1},\dots,\phi_{s}
\in {\cal{E}}_{G_{1}}(W)$. Assume that 
the ordered sequences $(a(x),b(x))$ and 
$(\psi_{1}(x),\dots,\psi_{r}(x), a(x), b(x), \phi_{1}(x), \dots, \phi_{s}(x))$ 
are compatible.
Let $k$ be a nonnegative integer such that
\begin{eqnarray}
& &(x_{1}-x_{2})^{k}\left(\prod_{1\le p<q\le r}(y_{p}-y_{q})^{k}\right)
\left(\prod_{i=1}^{r}(x_{1}-y_{i})^{k}(x_{2}-y_{i})^{k}\right)
\left(\prod_{j=1}^{s}(x_{1}-z_{j})^{k}(x_{2}-z_{j})^{k}\right)\nonumber\\
& &\cdot 
\left(\prod_{1\le i\le r, 1\le j\le s}(y_{i}-z_{j})^{k}\right)
\left(\prod_{1\le p<q\le s}(z_{p}-z_{q})^{k}\right)\nonumber\\
& &\;\;\cdot \psi_{1}(y_{1})\cdots\psi_{r}(y_{r})a(x_{1})b(x_{2})
\phi_{1}(z_{1})\cdots \phi_{s}(z_{s})\nonumber\\
&\in& \Hom (W,W((y_{1},\dots,y_{r},x_{1},x_{2},z_{1},\dots,z_{s}))).
\end{eqnarray}
Let $w\in W$ and let $l$ be a nonnegative integer such that
\begin{eqnarray}\label{epreparationequation}
&&x_{1}^{l}(x_{1}-x_{2})^{k}\left(\prod_{1\le p<q\le r}(y_{p}-y_{q})^{k}\right)
\left(\prod_{i=1}^{r}(x_{1}-y_{i})^{k}(x_{2}-y_{i})^{k}\right)
\left(\prod_{j=1}^{s}(x_{1}-z_{j})^{k}(x_{2}-z_{j})^{k}\right)\nonumber\\
&&\cdot \left(\prod_{1\le i\le r, 1\le j\le s}(y_{i}-z_{j})^{k}\right)
\left(\prod_{1\le p<q\le s}(z_{p}-z_{q})^{k}\right)\nonumber\\
&&\;\;\cdot\psi_{1}(y_{1})\cdots\psi_{r}(y_{r})a(x_{1})b(x_{2})
\phi_{1}(z_{1})\cdots \phi_{s}(z_{s})w\nonumber\\
&\in& W[[x_{1}]]((y_{1},\dots,y_{r},x_{2},z_{1},\dots,z_{s})).
\end{eqnarray}
Then
\begin{eqnarray}\label{epreparation1}
& &(x_{0}+x_{2})^{l}\left(\prod_{j=1}^{s}(x_{0}+x_{2}-z_{j})^{k}\right)
\psi_{1}(y_{1})\cdots\psi_{r}(y_{r})
(Y_{\cal{E}}(a,x_{0})b)(x_{2})\phi_{1}(z_{1})\cdots \phi_{s}(z_{s})w\nonumber\\
&=&(x_{0}+x_{2})^{l}\left(\prod_{j=1}^{s}(x_{0}+x_{2}-z_{j})^{k}\right)
\psi_{1}(y_{1})\cdots\psi_{r}(y_{r})a(x_{0}+x_{2})b(x_{2})\phi_{1}(z_{1})\cdots \phi_{s}(z_{s})w.
\end{eqnarray}
\el

{\pf} Set 
$$P=\prod_{1\le i<j\le r}(y_{i}-y_{j})^{k},\;\;\;\; 
Q=\prod_{1\le i<j\le s}(z_{i}-z_{j})^{k},\;\;\;\;
R=\prod_{1\le i\le r,\; 1\le j\le s}(y_{i}-z_{j})^{k}$$
and
$$S=\prod_{i=1}^{r}(x_{0}+x_{2}-y_{i})^{k}(x_{2}-y_{i})^{k}.$$
Let $m_{1},\dots,m_{s}$ be {\em arbitrarily fixed integers}. 
Since $\phi_{1}(z_{1})_{\ge m_{1}}\cdots \phi_{s}(z_{s})_{\ge m_{s}}w$ is a 
finite sum, from Proposition \ref{pdefinitionassociativity}
there exists a nonnegative integer $l'$ 
(depending on $m_{1},\dots,m_{s}$) such that
\begin{eqnarray}
& &(x_{0}+x_{2})^{l'}(Y_{\cal{E}}(a,x_{0})b)(x_{2})
\phi_{1}(z_{1})_{\ge m_{1}}\cdots \phi_{s}(z_{s})_{\ge m_{s}}w\nonumber\\
&=&(x_{0}+x_{2})^{l'}a(x_{0}+x_{2})b(x_{2})
\phi_{1}(z_{1})_{\ge m_{1}}\cdots \phi_{s}(z_{s})_{\ge m_{s}}w.
\end{eqnarray}
Then using (\ref{esimplefactresidule}) we get
\begin{eqnarray}
& &\left(\prod_{i=1}^{s}\Res_{z_{i}}z_{i}^{m_{i}}\right)
PQRS \left(\prod_{i=1}^{s}(x_{0}+x_{2}-z_{i})^{k}(x_{2}-z_{i})^{k}\right)
(x_{0}+x_{2})^{l'}\nonumber\\
& &\;\;\cdot \psi_{1}(y_{1})\cdots\psi_{r}(y_{r})(Y_{\cal{E}}(a,x_{0})b)(x_{2})\phi_{1}(z_{1})
\cdots \phi_{s}(z_{s})w\nonumber\\
&=&\left(\prod_{i=1}^{s}\Res_{z_{i}}z_{i}^{m_{i}}\right)PQRS
\left(\prod_{i=1}^{s}(x_{0}+x_{2}-z_{i})^{k}(x_{2}-z_{i})^{k}\right)\nonumber\\
& &\;\;\cdot (x_{0}+x_{2})^{l'}\psi_{1}(y_{1})\cdots\psi_{r}(y_{r})
(Y_{\cal{E}}(a,x_{0})b)(x_{2})\phi_{1}(z_{1})_{\ge m_{1}}
\cdots \phi_{s}(z_{s})_{\ge m_{s}}w\nonumber\\
&=&\left(\prod_{i=1}^{s}\Res_{z_{i}}z_{i}^{m_{i}}\right)PQRS
\left(\prod_{i=1}^{s}(x_{0}+x_{2}-z_{i})^{k}(x_{2}-z_{i})^{k}\right)
\nonumber\\
& &\;\;\cdot (x_{0}+x_{2})^{l'}\psi_{1}(y_{1})\cdots\psi_{r}(y_{r})
a(x_{0}+x_{2})b(x_{2})\phi_{1}(z_{1})_{\ge m_{1}}
\cdots \phi_{s}(z_{s})_{\ge m_{s}}w\nonumber\\
&=&\left(\prod_{i=1}^{s}\Res_{z_{i}}z_{i}^{m_{i}}\right)PQRS
\left(\prod_{i=1}^{s}(x_{0}+x_{2}-z_{i})^{k}(x_{2}-z_{i})^{k}\right)
\nonumber\\
& &\;\;\cdot (x_{0}+x_{2})^{l'}\psi_{1}(y_{1})\cdots\psi_{r}(y_{r})
a(x_{0}+x_{2})b(x_{2})\phi_{1}(z_{1})
\cdots \phi_{s}(z_{s})w.
\end{eqnarray}
Multiplying both sides by $(x_{0}+x_{2})^{l}$ we get
\begin{eqnarray}\label{eproofpreparation}
& &(x_{0}+x_{2})^{l'}\left(\prod_{i=1}^{s}\Res_{z_{i}}z_{i}^{m_{i}}\right)
PQRSx_{0}^{k}(x_{0}+x_{2})^{l}
\left(\prod_{i=1}^{s}(x_{0}+x_{2}-z_{i})^{k}(x_{2}-z_{i})^{k}\right)\nonumber\\
& &\cdot \psi_{1}(y_{1})\cdots\psi_{r}(y_{r})
(Y_{\cal{E}}(a,x_{0})b)(x_{2})\phi_{1}(z_{1})\cdots \phi_{s}(z_{s})w
\;\;\;\;\;\;\;\nonumber\\
&=&(x_{0}+x_{2})^{l'}\left(\prod_{i=1}^{s}\Res_{z_{i}}z_{i}^{m_{i}}\right)
 [PQRS x_{0}^{k}(x_{0}+x_{2})^{l}
\left(\prod_{i=1}^{s}(x_{0}+x_{2}-z_{i})^{k}(x_{2}-z_{i})^{k}\right)\nonumber\\
& &\cdot\psi_{1}(y_{1})\cdots\psi_{r}(y_{r})
a(x_{0}+x_{2})b(x_{2})\phi_{1}(z_{1})\cdots \phi_{s}(z_{s})w].
\end{eqnarray}
Noticing that from (\ref{epreparationequation}),
the expression in the bracket on the right-hand side of
(\ref{eproofpreparation}) involves only nonnegative 
powers of $x_{0}$, then multiplying (\ref{eproofpreparation})
by $x_{0}^{-k}(x_{2}+x_{0})^{-l'}$ we get
\begin{eqnarray}
& &\left(\prod_{i=1}^{s}\Res_{z_{i}}z_{i}^{m_{i}}\right) PQRS (x_{0}+x_{2})^{l}
\left(\prod_{i=1}^{s}(x_{0}+x_{2}-z_{i})^{k}(x_{2}-z_{i})^{k}\right)\nonumber\\
& &\;\;\cdot \psi_{1}(y_{1})\cdots\psi_{r}(y_{r})
(Y_{\cal{E}}(a,x_{0})b)(x_{2})\phi_{1}(z_{1})\cdots
\phi_{s}(z_{s})w\;\;\;\;\;\;\;\;\nonumber\\
&=&\left(\prod_{i=1}^{s}\Res_{z_{i}}z_{i}^{m_{i}}\right) PQRS (x_{0}+x_{2})^{l}
\left(\prod_{i=1}^{s}(x_{0}+x_{2}-z_{i})^{k}(x_{2}-z_{i})^{k}\right)
\nonumber\\
& &\;\;\cdot \psi_{1}(y_{1})\cdots\psi_{r}(y_{r})
a(x_{0}+x_{2})b(x_{2})\phi_{1}(z_{1})\cdots \phi_{s}(z_{s})w.
\end{eqnarray}
Since $l$ and $k$ are {\em independent of $m_{i}$'s and $m_{i}$'s are
arbitrary}, we have
\begin{eqnarray}
& &PQRS (x_{0}+x_{2})^{l}
\left(\prod_{i=1}^{s}(x_{0}+x_{2}-z_{i})^{k}(x_{2}-z_{i})^{k}\right)\nonumber\\
& &\cdot \psi_{1}(y_{1})\cdots\psi_{r}(y_{r})
(Y_{\cal{E}}(a,x_{0})b)(x_{2})\phi_{1}(z_{1})\cdots \phi_{s}(z_{s})w
\nonumber\\
&=&PQRS (x_{0}+x_{2})^{l}
\left(\prod_{i=1}^{s}(x_{0}+x_{2}-z_{i})^{k}(x_{2}-z_{i})^{k}\right)\nonumber\\
& &\cdot \psi_{1}(y_{1})\cdots\psi_{r}(y_{r})
a(x_{0}+x_{2})b(x_{2})\phi_{1}(z_{1})\cdots \phi_{s}(z_{s})w.\;\;\;\;\;\;
\end{eqnarray}
Noticing that we are allowed to multiply both sides by
$$\prod_{1\le p<q\le r}(y_{p}-y_{q})^{-k}
\prod_{1\le i\le r,1\le j\le s}(y_{i}-z_{j})^{-k}
\prod_{j=1}^{s}(x_{2}-z_{j})^{-k}\prod_{1\le p<q\le s} (z_{p}-z_{q})^{-k}$$
and by $\prod_{i=1}^{s}(-y_{i}+x_{0}+x_{2})^{k}(-y_{i}+x_{2})^{k}$
(but we are not allowed to multiply both sides by
$\prod_{i=1}^{s}(x_{0}+x_{2}-z_{i})^{-k}$),
we get the desired result.
$\;\;\;\;\Box$

Now we are in a position to prove our first key result:

\bt{tclosed}
Let $V$ be a subspace of ${\cal{E}}_{G_{1}}(W)$ such that any sequence
in $V$ of length $2$ or $3$ is compatible and such that
\begin{eqnarray}
& &1_{W}\in V,\\
& &\psi(x)_{n}\phi(x)\in V\;\;\;\mbox{ for }\psi(x),\phi(x)\in V,\; n\in \Z. 
\end{eqnarray}
Then $(V,Y_{\cal{E}},1_{W})$ carries the structure of a weak axiomatic 
$G_{1}$-vertex algebra with $W$ as a natural faithful module 
where the vertex operator map $Y_{W}$ is given by
$Y_{W}(\alpha(x),x_{0})=\alpha(x_{0})$.
Furthermore, assume that for any $\psi(x),\theta(x)\in V$, 
there exists a nonnegative integer
$k$ such that for every $\phi(x)\in V$ there exists 
a nonnegative integer $k'$ such that
\begin{eqnarray}
(x-y)^{k'}(y-z)^{k'}(x-z)^{k}\psi(x)\phi(y)\theta(z)
\in \Hom (W,W((x,y,z))).
\end{eqnarray}
Then $(V,Y_{\cal{E}},1_{W})$ carries the structure of 
an axiomatic $G_{1}$-vertex algebra. 
\et

{\pf} For the assertion on the axiomatic $G_{1}$-vertex algebra structure,
with Proposition \ref{pnthproduct} and
Lemma \ref{lvacuum-creationproperties}
we must prove the weak associativity, i.e.,
for $\psi,\phi,\theta\in V$, there exists a nonnegative 
integer $k$ such that
\begin{eqnarray}\label{eweakassocmainthem}
(x_{0}+x_{2})^{k}Y_{\cal{E}}(\psi,x_{0}+x_{2})Y_{\cal{E}}(\phi,x_{2})\theta
=(x_{0}+x_{2})^{k}Y_{\cal{E}}(Y_{\cal{E}}(\psi,x_{0})\phi, x_{2})\theta.
\end{eqnarray}
Let $k$ and $k'$ be nonnegative integers such that
\begin{eqnarray}
(x-y)^{k'}(x-z)^{k}(y-z)^{k'}
\psi(x)\phi(y)\theta(z)\in \Hom (W,W((x,y,z))).
\end{eqnarray}
For the first assertion, both $k$ and $k'$ depend on all $\psi,\phi,\theta$
and for the second assertion, $k$ depends only on $\psi$ and $\theta$.

Let $w\in W$ be {\em arbitrary and fixed}. 
There exists a nonnegative integer $l$ such that
\begin{eqnarray}\label{elllkkkw}
x^{l}y^{l}z^{l}(x-y)^{k'}(x-z)^{k}(y-z)^{k'}
\psi(x)\phi(y)\theta(z)w\in W[[x,y,z]].
\end{eqnarray}
In view of Proposition \ref{pdefinitionassociativity},
by replacing $l$ with a larger integer if necessary we may assume that
\begin{eqnarray}
(x_{2}+x)^{l}\phi(x_{2}+x)\theta(x)w
=(x_{2}+x)^{l}(Y_{\cal{E}}(\phi,x_{2})\theta)(x)w.
\end{eqnarray}
Then
\begin{eqnarray}
& &x^{l}x_{0}^{k}(x_{0}+x)^{l}x_{2}^{k'}(x_{0}-x_{2})^{k'}(x_{2}+x)^{l}
\psi(x_{0}+x)\phi(x_{2}+x)\theta(x)w\nonumber\\
&=&x^{l}x_{0}^{k}(x_{0}+x)^{l}x_{2}^{k'}(x_{0}-x_{2})^{k'}(x_{2}+x)^{l}
\psi(x_{0}+x)(Y_{\cal{E}}(\phi,x_{2})\theta)(x)w.
\end{eqnarray}
Noticing that the expression on the left-hand side lies in $W[[x,x_{0},x_{2}]]$
by (\ref{elllkkkw}), we have
\begin{eqnarray}
x^{l}x_{0}^{k}(x_{0}+x)^{l}x_{2}^{k'}(x_{0}-x_{2})^{k'}(x_{2}+x)^{l}
\psi(x_{0}+x)(Y_{\cal{E}}(\phi,x_{2})\theta)(x)w
\in W[[x,x_{0},x_{2}]].
\end{eqnarray}
By multiplying by 
$x^{-l}x_{0}^{-k}x_{2}^{-k'}(x_{0}-x_{2})^{-k'}(x+x_{2})^{-l}$, which lies in 
$\C ((x,x_{0}))((x_{2}))$,
we have
\begin{eqnarray}
(x_{0}+x)^{l}
\psi(x_{0}+x)(Y_{\cal{E}}(\phi,x_{2})\theta)(x)w
\in W((x,x_{0}))((x_{2})).
\end{eqnarray}
In view of Proposition \ref{pdefinitionassociativity}, 
by considering components of $Y_{\cal{E}}(\phi,x_{2})\theta$,
we have
\begin{eqnarray}
(x_{0}+x)^{l}\psi(x_{0}+x)(Y_{\cal{E}}(\phi,x_{2})\theta)(x)w
=(x_{0}+x)^{l}
(Y_{\cal{E}}(\psi,x_{0})Y_{\cal{E}}(\phi,x_{2})\theta)(x)w.
\end{eqnarray}
Then
\begin{eqnarray}\label{e5.76}
& &(x_{0}+x_{2})^{k}(x_{2}+x)^{l}(x_{0}+x_{2}+x)^{l}
\psi(x_{0}+x_{2}+x)\phi(x_{2}+x)\theta(x)w\nonumber\\
&=&(x_{0}+x_{2})^{k}(x_{2}+x)^{l}(x_{0}+x_{2}+x)^{l}
\psi(x_{0}+x_{2}+x)(Y_{\cal{E}}(\phi,x_{2})\theta)(x)w\nonumber\\
&=&e^{x_{2}{\partial\over\partial x_{0}}}x_{0}^{k}(x_{2}+x)^{l}(x_{0}+x)^{l}
\psi(x_{0}+x)(Y_{\cal{E}}(\phi,x_{2})\theta)(x)w\nonumber\\
&=&e^{x_{2}{\partial\over\partial x_{0}}}x_{0}^{k}(x_{2}+x)^{l}(x_{0}+x)^{l}
(Y_{\cal{E}}(\psi,x_{0})Y_{\cal{E}}(\phi,x_{2})\theta)(x)w\nonumber\\
&=&(x_{0}+x_{2})^{k}(x_{2}+x)^{l}(x_{0}+x_{2}+x)^{l}
(Y_{\cal{E}}(\psi,x_{0}+x_{2})Y_{\cal{E}}(\phi,x_{2})\theta)(x)w.
\end{eqnarray}

On the other hand, let $n\in \Z$ be {\em arbitrarily fixed}.
Since $\psi(x)_{m}\phi(x)=0$ for $m$ sufficiently large, there exists
a nonnegative integer $l'\in \Z$ such that
\begin{eqnarray}\label{eges}
(x_{2}+x)^{l'}(Y_{\cal{E}}(\psi(x)_{m}\phi(x), x_{2})\theta(x))w
=(x_{2}+x)^{l'}(\psi(x)_{m}\phi(x))(x_{2}+x)\theta(x)w
\end{eqnarray}
for {\em all} $m\ge n$. With (\ref{elllkkkw}),
in view of Lemma \ref{lbeforegenerating}, we have
\begin{eqnarray}\label{ethis}
& &(x_{0}+x_{2})^{l}(x_{0}+x_{2}-x)^{k}\psi(x_{0}+x_{2})\phi(x_{2})\theta(x)w
\nonumber\\
&=&(x_{0}+x_{2})^{l}(x_{0}+x_{2}-x)^{k}
(Y_{\cal{E}}(\psi,x_{0})\phi)(x_{2})\theta(x)w.
\end{eqnarray}
Using (\ref{esimplefactresidule}), (\ref{eges}) and (\ref{ethis}) we get
\begin{eqnarray}\label{e5.78}
& &\Res_{x_{0}}x_{0}^{n}(x_{0}+x_{2}+x)^{l}(x_{0}+x_{2})^{k}
(x_{2}+x)^{l'}
(Y_{\cal{E}}(Y_{\cal{E}}(\psi,x_{0})\phi, x_{2})\theta)(x)w\nonumber\\
&=&\Res_{x_{0}}x_{0}^{n}(x_{0}+x_{2}+x)^{l}(x_{0}+x_{2})^{k}
(x_{2}+x)^{l'}
(Y_{\cal{E}}(Y_{\cal{E}}(\psi,x_{0})_{\ge n}\phi, x_{2})\theta)(x)w\nonumber\\
&=&\Res_{x_{0}}x_{0}^{n}(x_{0}+x_{2}+x)^{l}
(x_{0}+x_{2})^{k}(x_{2}+x)^{l'}
(Y_{\cal{E}}(\psi,x_{0})_{\ge n}\phi)(x_{2}+x)\theta(x)w\nonumber\\
&=&\Res_{x_{0}}x_{0}^{n}(x_{0}+x_{2}+x)^{l}
(x_{0}+x_{2})^{k}(x_{2}+x)^{l'}
(Y_{\cal{E}}(\psi,x_{0})\phi)(x_{2}+x)\theta(x)w\nonumber\\
&=&\Res_{x_{0}}x_{0}^{n}(x_{2}+x)^{l'}e^{x{\partial\over\partial x_{2}}}
\left[(x_{0}+x_{2})^{l}(x_{0}+x_{2}-x)^{k}
(Y_{\cal{E}}(\psi,x_{0})\phi)(x_{2})\theta(x)w\right]\nonumber\\
&=&\Res_{x_{0}}x_{0}^{n}(x_{2}+x)^{l'}e^{x{\partial\over\partial x_{2}}}
\left[(x_{0}+x_{2})^{l}(x_{0}+x_{2}-x)^{k}
\psi(x_{0}+x_{2})\phi(x_{2})\theta(x)w\right]\nonumber\\
&=&\Res_{x_{0}}x_{0}^{n}(x_{2}+x)^{l'}(x_{0}+x_{2}+x)^{l}(x_{0}+x_{2})^{k}
\psi(x_{0}+x_{2}+x)\phi(x_{2}+x)\theta(x)w.\;\;\;\;\;\;\;\;
\end{eqnarray}
Combining (\ref{e5.78}) with (\ref{e5.76}) we get
\begin{eqnarray}\label{e5.79}
& &\Res_{x_{0}}x_{0}^{n}(x_{0}+x_{2}+x)^{l}
(x_{0}+x_{2})^{k}(x_{2}+x)^{l+l'}
(Y_{\cal{E}}(\psi,x_{0}+x_{2})Y_{\cal{E}}(\phi,x_{2})\theta)(x)w\nonumber\\
&=&\Res_{x_{0}}x_{0}^{n}(x_{0}+x_{2}+x)^{l}(x_{0}+x_{2})^{k}
(x_{2}+x)^{l+l'}
(Y_{\cal{E}}(Y_{\cal{E}}(\psi,x_{0})\phi, x_{2})\theta)(x)w.
\end{eqnarray}
Notice that both sides of (\ref{e5.79}) involve only finitely 
many negative powers of $x_{2}$.
Then multiplying both sides by $(x+x_{2})^{-l-l'}$ we get
\begin{eqnarray}
& &\Res_{x_{0}}x_{0}^{n}
(x_{0}+x_{2}+x)^{l}(x_{0}+x_{2})^{k}
(Y_{\cal{E}}(\psi,x_{0}+x_{2})Y_{\cal{E}}(\phi,x_{2})\theta)(x)w\nonumber\\ 
&=&\Res_{x_{0}}x_{0}^{n}(x_{0}+x_{2}+x)^{l}(x_{0}+x_{2})^{k}
(Y_{\cal{E}}(Y_{\cal{E}}(\psi,x_{0})\phi, x_{2})\theta)(x)w.
\end{eqnarray}
Since {\em $n$ is arbitrary and $l$ and $k$ do not depend on $n$}, we must have
\begin{eqnarray}\label{enearfinal}
& &(x_{0}+x_{2}+x)^{\ell}(x_{0}+x_{2})^{k}
(Y_{\cal{E}}(\psi,x_{0}+x_{2})Y_{\cal{E}}(\phi,x_{2})\theta)(x)w\nonumber\\
&=&(x_{0}+x_{2}+x)^{\ell}(x_{0}+x_{2})^{k}
(Y_{\cal{E}}(Y_{\cal{E}}(\psi,x_{0})\phi, x_{2})\theta)(x)w.
\end{eqnarray}
Notice that $(x+x_{0}+x_{2})^{-l}=(x+x_{2}+x_{0})^{-l}$ and that
we are allowed to multiply the left-hand side of (\ref{enearfinal})
by $(x+x_{0}+x_{2})^{-l}$
and to multiply the right-hand side by $(x+x_{2}+x_{0})^{-l}$.
Then multiplying both sides by $(x+x_{0}+x_{2})^{-l}$ 
we obtain
\begin{eqnarray}\label{e4.39}
& &(x_{0}+x_{2})^{k}(Y_{\cal{E}}(\psi,x_{0}+x_{2})Y_{\cal{E}}(\phi,x_{2})
\theta)(x)w\nonumber\\
&=&(x_{0}+x_{2})^{k}
(Y_{\cal{E}}(Y_{\cal{E}}(\psi,x_{0})\phi,x_{2})\theta)(x)w.
\end{eqnarray}
Since $k$ does not depend on $w$, we immediately have
(\ref{eweakassocmainthem}), as desired.

For $a(x),b(x)\in V,\; w\in W$, in view of 
Proposition \ref{pdefinitionassociativity} there exists a
nonnegative integer $l$ such that
$$(x_{0}+x_{2})^{l}(Y_{\cal{E}}(a(x),x_{0})b(x))(x_{2})w=
(x_{0}+x_{2})^{l}a(x_{0}+x_{2})b(x_{2})w.$$
That is,
\begin{eqnarray}
(x_{0}+x_{2})^{l}Y_{W}(Y_{\cal{E}}(a(x),x_{0})b(x),x_{2})w
=(x_{0}+x_{2})^{l}Y_{W}(a(x),x_{0}+x_{2})Y_{W}(b(x),x_{2})w.
\end{eqnarray}
Therefore $W$ is a $V$-module with $Y_{W}(\alpha(x),x_{0})=\alpha(x_{0})$
for $\alpha(x)\in V$.
$\;\;\;\;\Box$

Our next goal is to prove that any compatible set $S$ of 
$G_{1}$-vertex operators on $W$ gives rise to an axiomatic
$G_{1}$-vertex algebra.
To achieve this goal, we first need to show that 
for $a,b,c\in S,\; n\in \Z$, the sequences $(a(x)_{n}b(x),c(x))$
and $(c(x),a(x)_{n}b(x))$ are compatible,
so that $c(x)_{m}(a(x)_{n}b(x))$ and
$(a(x)_{n}b(x))_{m}c(x)$ are defined for $m\in \Z$. 
The following is another key result:

\bp{pgeneratingcomplicatedone}
Let $\psi_{1}(x),\dots,\psi_{r}(x), a(x),b(x),\phi_{1}(x), \dots,\phi_{s}(x)
\in {\cal{E}}_{G_{1}}(W)$. Assume that 
the ordered sequences $(a(x), b(x))$ and 
$$(\psi_{1}(x),\dots,\psi_{r}(x), a(x),b(x),\phi_{1}(x), \dots,\phi_{s}(x))$$
 are compatible.
Then for any $n\in \Z$, the ordered sequence 
$$(\psi_{1}(x),\dots,\psi_{r}(x),a(x)_{n}b(x),\phi_{1}(x),\dots,\phi_{s}(x))$$
 is compatible.
\ep

{\bf Proof.} Let $n\in \Z$ be {\em arbitrarily fixed}. 
{}From Proposition \ref{pdefinitionassociativity} 
there exists a nonnegative integer $k'$ such that 
\begin{eqnarray}\label{etruncationpsiphi}
x_{0}^{k'+n}Y_{\cal{E}}(a,x_{0})b\in {\cal{E}}_{G_{1}}(W)[[x_{0}]].
\end{eqnarray}
Let $k$ be a nonnegative integer such that
\begin{eqnarray}
& &\left(\prod_{1\le i<j\le r}(y_{i}-y_{j})^{k}\right)
\left(\prod_{1\le i\le r, 1\le j\le s}(y_{i}-z_{j})^{k}\right)
\left(\prod_{1\le i<j\le s}(z_{i}-z_{j})^{k}\right)\nonumber\\
& &\;\;\cdot (x_{1}-x_{2})^{k}
\left(\prod_{i=1}^{r}(x_{1}-y_{i})^{k}(x_{2}-y_{i})^{k}\right)
\left(\prod_{i=1}^{s}(x_{1}-z_{i})^{k}(x_{2}-z_{i})^{k}\right)\nonumber\\
& &\;\;\cdot \psi_{1}(y_{1})\cdots \psi_{r}(y_{r})
a(x_{1})b(x_{2})\phi_{1}(z_{1})\cdots \phi_{s}(z_{s})\nonumber\\
& &\in \Hom (W,W((y_{1},\dots, y_{r},x_{1},x_{2},z_{1},\dots,z_{s}))).
\end{eqnarray}
Set
$$P=\prod_{1\le i<j\le r}(y_{i}-y_{j})^{k},\;\;\;\; 
Q=\prod_{1\le i<j\le s}(z_{i}-z_{j})^{k},\;\;\;\;
R=\prod_{1\le i\le r,\; 1\le j\le s}(y_{i}-z_{j})^{k}.$$
Let $w\in W$. 
Then there exists a nonnegative integer $l$ such that
\begin{eqnarray}
& &x_{1}^{l} P Q R (x_{1}-x_{2})^{k}
\left(\prod_{i=1}^{r}(x_{1}-y_{i})^{k}(x_{2}-y_{i})^{k}\right)
\left(\prod_{j=1}^{s}(x_{1}-z_{j})^{k}(x_{2}-z_{j})^{k}\right)\nonumber\\
& &\;\;\cdot\psi_{1}(y_{1})\cdots \psi_{r}(y_{r})a(x_{1})b(x_{2})
\phi_{1}(z_{1})\cdots \phi_{s}(z_{s})w\nonumber\\
&\in& W[[x_{1}]]((x_{2},y_{1},\dots, y_{r},z_{1},\dots,z_{s}))
\end{eqnarray}
(cf. Lemma \ref{lbeforegenerating}). Hence
\begin{eqnarray}
& &PQR(x_{0}+x_{2})^{l}x_{0}^{k}
\left(\prod_{i=1}^{r}(x_{0}+x_{2}-y_{i})^{k}(x_{2}-y_{i})^{k}\right)
\left(\prod_{j=1}^{s}(x_{0}+x_{2}-z_{j})^{k}(x_{2}-z_{j})^{k}\right)
\nonumber\\
& &\;\;\cdot\psi_{1}(y_{1})\cdots \psi_{r}(y_{r})a(x_{0}+x_{2})b(x_{2})
\phi_{1}(z_{1})\cdots \phi_{s}(z_{s})w\nonumber\\
&\in& W[[x_{0}]]((x_{2},y_{1},\dots, y_{r},z_{1},\dots,z_{s})).
\end{eqnarray}
In the following we are going to use the binomial expansions 
for $x_{2}^{l+k'}=((x_{2}+x_{0})-x_{0})^{l+k'}$,
$(x_{2}-y_{i})^{k+k'}=((x_{2}+x_{0}-y_{i})-x_{0})^{k+k'}$ and
$(x_{2}-z_{j})^{k+k'}=((x_{2}+x_{0}-z_{j})-x_{0})^{k+k'}$.
Using (\ref{etruncationpsiphi}) and Lemma \ref{lbeforegenerating}
we obtain
\begin{eqnarray}\label{ecompatibilitythreeproof}
& &x_{2}^{l+k'}\prod_{i=1}^{r}(x_{2}-y_{i})^{2k+k'}
\prod_{j=1}^{s}(x_{2}-z_{i})^{2k+k'}\nonumber\\
& &\cdot \psi_{1}(y_{1})\cdots\psi_{r}(y_{r})
(a(x)_{n}b(x))(x_{2})\phi_{1}(z_{1})\cdots\phi_{s}(z_{s})w\nonumber\\
&=&\Res_{x_{0}}x_{0}^{n}x_{2}^{l+k'}\prod_{i=1}^{r}(x_{2}-y_{i})^{2k+k'}
\prod_{j=1}^{s}(x_{2}-z_{i})^{2k+k'}\nonumber\\
& &\cdot \psi_{1}(y_{1})\cdots \psi_{r}(y_{r})
(Y_{\cal{E}}(a,x_{0})b)(x_{2})\phi_{1}(z_{1})\cdots \phi_{s}(z_{s})w\nonumber\\
&=&\Res_{x_{0}}x_{0}^{n}\left(\sum_{p\ge 0}{l+k'\choose p}
(x_{0}+x_{2})^{l+k'-p}(-x_{0})^{p}\right)\nonumber\\
& &\cdot \prod_{i=1}^{r}\left(\sum_{t\ge 0}{k+k'\choose t}(x_{2}+x_{0}-y_{i})^{k+k'-t}
(-x_{0})^{t}\right)\left(\prod_{i=1}^{r}(x_{2}-y_{i})^{k}\right)
\nonumber\\
& &\cdot \prod_{j=1}^{s}\left(\sum_{q\ge 0}{k+k'\choose q}
(x_{0}+x_{2}-z_{j})^{k+k'-q}(-x_{0})^{q}\right)
\left(\prod_{j=1}^{s}(x_{2}-z_{j})^{k}\right)\nonumber\\
& &\cdot \psi_{1}(y_{1})\cdots \psi_{r}(y_{r})
(Y_{\cal{E}}(a,x_{0})b)(x_{2})\phi_{1}(z_{1})\cdots \phi_{s}(z_{s})w\nonumber\\
&=&\Res_{x_{0}}x_{0}^{n}\left(\sum_{0\le p\le k'}{l+k'\choose p}
(x_{0}+x_{2})^{l+k'-p}(-x_{0})^{p}\right)\nonumber\\
& &\cdot \prod_{i=1}^{r}\left(\sum_{0\le t\le k'}{k+k'\choose t}(x_{2}+x_{0}-y_{i})^{k+k'-t}
(-x_{0})^{t}\right)\left(\prod_{i=1}^{r}(x_{2}-y_{i})^{k}\right)\nonumber\\
& &\cdot \prod_{j=1}^{s}\left(\sum_{0\le q\le k'}{k+k'\choose q}
(x_{0}+x_{2}-z_{j})^{k+k'-q}(-x_{0})^{q}\right)
\left(\prod_{j=1}^{s}(x_{2}-z_{j})^{k}\right)\nonumber\\
& &\cdot \psi_{1}(y_{1})\cdots \psi_{r}(y_{r})(Y_{\cal{E}}(a,x_{0})b)(x_{2})
\phi_{1}(z_{1})\cdots \phi_{s}(z_{s})w\nonumber\\
&=&\Res_{x_{0}}x_{0}^{n}\left(\sum_{0\le p\le k'}{l+k'\choose p}
(x_{0}+x_{2})^{l+k'-p}(-x_{0})^{p}\right)\nonumber\\
& &\cdot \prod_{i=1}^{r}\left(\sum_{0\le t\le k'}{k+k'\choose t}(x_{2}+x_{0}-y_{i})^{k+k'-t}
(-x_{0})^{t}\right)\left(\prod_{i=1}^{r}(x_{2}-y_{i})^{k}\right)\nonumber\\
& &\cdot \prod_{j=1}^{s}\left(\sum_{0\le q\le k'}{k+k'\choose q}
(x_{0}+x_{2}-z_{j})^{k+k'-q}(-x_{0})^{q}\right)
\left(\prod_{j=1}^{s}(x_{2}-z_{j})^{k}\right)\nonumber\\
& &\cdot \psi_{1}(y_{1})\cdots \psi_{r}(y_{r})
a(x_{0}+x_{2})b(x_{2})\phi_{1}(z_{1})\cdots \phi_{s}(z_{s})w.
\end{eqnarray}
Noticing that
\begin{eqnarray}
& &PQR \left(\prod_{i=1}^{r}(x_{2}-y_{i})^{k}\right)
\left(\prod_{j=1}^{s}(x_{2}-z_{j})^{k}\right)
(x_{0}+x_{2})^{l+k'-p}
\left(\prod_{i=1}^{r}(x_{2}+x_{0}-y_{i})^{k+k'-t}\right)
\nonumber\\
& &\cdot \left(\prod_{j=1}^{s}(x_{0}+x_{2}-z_{j})^{k+k'-q}\right)
\psi_{1}(y_{1})\cdots \psi_{r}(y_{r})
a(x_{0}+x_{2})b(x_{2})\phi_{1}(z_{1})\cdots \phi_{s}(z_{s})w\nonumber\\
&\in& W((y_{1},\dots,y_{r},x_{0},x_{2},z_{1},\dots,z_{s}))
\end{eqnarray}
for $0\le p,t, q\le k'$, from (\ref{ecompatibilitythreeproof})  we have that
\begin{eqnarray*}
& &P Q R x_{2}^{l+k'}\prod_{i=1}^{r}(x_{2}-y_{i})^{2k+k'}
\prod_{j=1}^{s}(x_{2}-z_{j})^{2k+k'}\nonumber\\
& &\;\;\cdot\psi_{1}(y_{1})\cdots \psi_{r}(y_{r})
(a(x)_{n}b(x))(x_{2})\phi_{1}(z_{1})\cdots \phi_{s}(z_{s})w
\end{eqnarray*}
lies in $W((y_{1},\dots, y_{r},x_{2},z_{1},\dots,z_{s}))$, and so
\begin{eqnarray}
& &P Q R \prod_{i=1}^{r}(x_{2}-y_{i})^{2k+k'}
\prod_{j=1}^{s}(x_{2}-z_{j})^{2k+k'}\nonumber\\
& &\;\;\cdot \psi_{1}(y_{1})\cdots \psi_{r}(y_{r})
(a(x)_{n}b(x))(x_{2})\phi_{1}(z_{1})\cdots \phi_{s}(z_{s})w
\;\;\;\;\;\;\nonumber\\
&\in& W((y_{1},\dots, y_{r},x_{2},z_{1},\dots,z_{s})).
\end{eqnarray}
This proves that the sequence 
$(\psi_{1}(x),\dots,\psi_{r}(x),a(x)_{n}b(x), \phi_{1}(x),
\dots,\phi_{s}(x))$ is compatible, since $k$ and $k'$ are independent of $w$.
$\;\;\;\;\Box$

It follows from Proposition \ref{pgeneratingcomplicatedone}
that if $S$ is a compatible set of weak $G_{1}$-vertex operators
on $W$, then for any
$a^{(1)},\dots,a^{(r+1)}\in S$, the expression
$$Y_{\cal{E}}(a^{(1)}(x),x_{1})\cdots 
Y_{\cal{E}}(a^{(r)}(x),x_{r})a^{(r+1)}(x)$$
is recursively well defined. 

The following result states that any maximal compatible subspace 
of ${\cal{E}}_{G_{1}}(W)$ is automatically closed and it has a
weak axiomatic $G_{1}$-vertex algebra structure.

\bp{pmaximal}
Let $V$ be a maximal compatible subspace of ${\cal{E}}_{G_{1}}(W)$.
Then $1_{W}\in V$ and 
\begin{eqnarray}
a(x)_{n}b(x)\in V\;\;\;\mbox{ for }a(x),b(x)\in V,\; n\in \Z.
\end{eqnarray}
Furthermore,
$(V,Y_{\cal{E}},1_{W})$ carries the structure of a weak axiomatic
$G_{1}$-vertex algebra with $W$ as a natural module where the vertex
operator map $Y_{W}$ is given by $Y_{W}(\alpha(x),x_{0})=\alpha(x_{0})$.
\ep

{\bf Proof.} Clearly the space spanned by $V$ and 
$1_{W}$ is still compatible. With $V$ being
maximal we must have $1_{W}\in V$.
Now, let $a(x),b(x)\in V$ and $n\in \Z$. In view of 
Proposition \ref{pgeneratingcomplicatedone},
any (ordered) sequence in $V\cup \{ a(x)_{n}b(x)\}$
with one appearance of $a(x)_{n}b(x)$ is compatible.
It follows from induction on the number of 
appearance of $a(x)_{n}b(x)$ and from
Proposition \ref{pgeneratingcomplicatedone} that
any (ordered) sequence in $V\cup \{ a(x)_{n}b(x)\}$
with any (finite) number of appearance of $a(x)_{n}b(x)$ is compatible.
So the space spanned by $V$ and $a(x)_{n}b(x)$ is compatible.
Again, with $V$ being maximal we must have $a(x)_{n}b(x)\in V$. 
This proves that $V$ is closed, and hence by
Theorem \ref{tclosed} $(V,Y_{\cal{E}},1_{W})$ carries the structure 
of a weak axiomatic $G_{1}$-vertex algebra with $W$ as a natural module.
$\;\;\;\;\Box$

Let $S$ be a compatible set of $G_{1}$-vertex operators on $W$.
By Zorn's lemma there exists a maximal compatible space $V$ of 
${\cal{E}}_{G_{1}}(W)$, containing $S$ and $1_{W}$, and then
by Proposition \ref{pmaximal}
$(V,Y_{\cal{E}},1_{W})$ carries the structure of a weak axiomatic
$G_{1}$-vertex algebra with $W$ as a natural module.
Now, $S$ as a subset of $V$ generates a subalgebra $\<S\>$
of $V$. Then in view of Proposition 2.22 we obtain 
our main result (cf. [B2], Theorem 7.9):

\bt{tgeneratingthem}
Let $S$ be any compatible set of $G_{1}$-vertex operators on $W$.
Then for any $a^{(1)}(x),\dots,a^{(r)}(x)\in S$, the expression 
$$Y_{\cal{E}}(a^{(1)}(x),x_{1})\cdots Y_{\cal{E}}(a^{(r)}(x),x_{r})1_{W}$$
is recursively well defined.
Furthermore, if we set 
\begin{eqnarray}
U={\rm span }\{ a^{(1)}(x)_{n_{1}}\cdots a^{(r)}(x)_{n_{r}}1_{W}
\;|\; r\ge 0,\; a^{(i)}(x)\in S,\; n_{i}\in \Z\},
\end{eqnarray}
then $(U,Y_{\cal{E}},1_{W})$
carries the structure of a weak
axiomatic $G_{1}$-vertex algebra with $W$ as a natural module
where the vertex
operator map $Y_{W}$ is given by $Y_{W}(\alpha(x),x_{0})=\alpha(x_{0})$.
$\;\;\;\;\Box$
\et

As an immediate consequence of Proposition \ref{ppracticalcase}
and Theorem \ref{tgeneratingthem} we have:

\bc{cgeneratingthemquasi}
Let $S$ be a set of weak $G_{1}$-vertex operators on $W$ such that
for any $a,b\in S$, there exists a nonnegative integer
$k$ such that 
\begin{eqnarray}
(x_{1}-x_{2})^{k}a(x_{1})b(x_{2})
=\sum_{i=1}^{r}\alpha_{i}(x_{1}-x_{2})^{k}b^{(i)}(x_{2})a^{(i)}(x_{1})
\end{eqnarray}
for some $\alpha_{i}\in \C,\; a^{(i)}, b^{(i)}\in S,\; r\ge 1$.
Then all the assertions of Theorem \ref{tgeneratingthem} hold.
$\;\;\;\;\Box$
\ec

Recall from Corollary \ref{clocalimplycompatibility}
that any space of pairwise mutually local vertex operators on
$W$ is compatible. Then in view of Theorem \ref{tclosed},
any closed space of pairwise mutually local vertex operators on
$W$ is a weak axiomatic $G_{1}$-vertex algebra with $W$ as a module.
In [Li1], it was proved (cf. [MN]) that any closed space of 
pairwise mutually local vertex operators on $W$ is 
an (ordinary) vertex algebra with $W$ as a module. 
Theorem \ref{tclosed} does not directly
imply the corresponding result of [Li1], but
Theorem \ref{tclosed} together with 
Propositions \ref{pcommutativesub} and \ref{pcommutant} does.

\bt{tgeneratingthemlocal}
Let $S$ be a set of pairwise mutually local
$G_{1}$-vertex operators on $W$. 
Then for $r\ge 0,\; a^{(i)}\in S$, $i=1,\dots,r$,
the expression 
$$Y_{\cal{E}}(a^{(1)}(x),x_{1})\cdots Y_{\cal{E}}(a^{(r)}(x),x_{r})1_{W}$$
is recursively well defined and
\begin{eqnarray}\label{eoldnew}
Y_{\cal{E}}(a^{(1)}(x),x_{1})\cdots Y_{\cal{E}}(a^{(r)}(x),x_{r})1_{W}
=Y_{\cal{E}}^{c}(a^{(1)}(x),x_{1})\cdots 
Y_{\cal{E}}^{c}(a^{(r)}(x),x_{r})1_{W},
\end{eqnarray}
where
\begin{eqnarray}
Y_{\cal{E}}^{c}(\alpha(x),x_{0})\beta(x)
=\Res_{x_{1}}\left(x_{0}^{-1}\delta\left(\frac{x_{1}-x}{x_{0}}\right)
\alpha(x_{1})\beta(x)-
x_{0}^{-1}\delta\left(\frac{x-x_{1}}{-x_{0}}\right)\beta(x)\alpha(x_{1})\right)
\end{eqnarray}
for $\alpha(x),\beta(x)\in {\cal{E}}_{G_{1}}(W)$.
Furthermore, if we set
\begin{eqnarray}
U={\rm span }\{ a^{(1)}(x)_{n_{1}}\cdots a^{(r)}(x)_{n_{r}}1_{W}
\;|\; r\ge 0,\; a^{(i)}\in S,\; n_{i}\in \Z\},
\end{eqnarray}
then $(U,Y_{\cal{E}}^{c},1_{W})$ carries the structure of 
an ordinary vertex algebra with $W$ as a natural module
where the vertex operator map $Y_{W}$ is given by 
$Y_{W}(\alpha(x),x_{0})=\alpha(x_{0})$.
\et

{\bf Proof.} Let $U'$ be the space defined in Theorem \ref{tgeneratingthem}.
In view of Corollary \ref{cgeneratingthemquasi}, $U'$ is a weak 
axiomatic $G_{1}$-vertex algebra with $W$ as a natural module.
Since $W$ is a faithful module and for any $a(x),b(x)\in S$,
$Y_{W}(a(x),x_{1})$ $(=a(x_{1}))$ and $Y_{W}(b(x),x_{2})$ $(=b(x_{2}))$
are mutually local, by Proposition \ref{pcommutant}, 
$Y(a(x),x_{1})$ and $Y(b(x),x_{2})$ acting on $U'$
are mutually local. By Proposition \ref{pcommutativesub},
$U'$ is an ordinary vertex algebra because $S$ generates $U'$ as
a weak axiomatic $G_{1}$-vertex algebra. In view of 
Proposition \ref{prestrictedgvamodule},
$W$ equipped with the linear map
$Y_{W}$ given by $Y_{W}(\alpha(x),x_{2})=\alpha(x_{2})$
is a module for $U'$ viewed as a vertex algebra.
Then for $\alpha(x), \beta(x)\in U'$,
\begin{eqnarray}
& &Y_{\cal{E}}(\alpha(x_{2}),x_{0})\beta(x_{2})\nonumber\\
&=&Y_{W}(Y_{\cal{E}}(\alpha,x_{0})\beta,x_{2})\nonumber\\
&=&\Res_{x_{1}}x_{0}^{-1}\delta\left(\frac{x_{1}-x_{2}}{x_{0}}\right)
Y_{W}(\alpha(x),x_{1})Y_{W}(\beta(x),x_{2})\nonumber\\
& &-\Res_{x_{1}}
x_{0}^{-1}\delta\left(\frac{x_{2}-x_{1}}{-x_{0}}\right)
Y_{W}(\beta(x),x_{2})Y_{W}(\alpha(x),x_{1})\nonumber\\
&=&\Res_{x_{1}}\left(x_{0}^{-1}\delta\left(\frac{x_{1}-x_{2}}{x_{0}}\right)
\alpha(x_{1})\beta(x_{2})-
x_{0}^{-1}\delta\left(\frac{x_{2}-x_{1}}{-x_{0}}\right)
\beta(x_{2})\alpha(x_{1})\right)
\nonumber\\
&=&Y_{\cal{E}}^{c}(\alpha(x_{2}),x_{0})\beta(x_{2}).
\end{eqnarray}
Now, (\ref{eoldnew}) follows immediately from induction and we have $U=U'$.
This completes the proof.$\;\;\;\;\Box$


\begin{thebibliography}{DLinM}

\bibitem[B1]{b1}
R. E. Borcherds, Vertex algebras, Kac-Moody algebras, and the Monster,
{\em Proc. Natl. Acad. Sci. USA} {\bf 83} (1986), 3068-3071.

\bibitem[B2]{b2}
R. E. Borcherds, Vertex algebras, in ``Topological Field Theory,
Primitive Forms and Related Topics'' (Kyoto, 1996), edited by
M. Kashiwara, A. Matsuo, K. Saito and I. Satake, Progress in Math.,
160, Birkh\"auser, Boston, 1998, 35-77; q-alg/9706008.

\bibitem[BPZ]{bpz}
A. A. Belavin, A. M. Polyakov and A. B. Zamolodchikov, 
Infinite conformal symmetries in two-dimensional quantum 
field theory, {\em Nucl. Phys.} {\bf B241} (1984), 333-380.

\bibitem[DL]{dl}
C. Dong and J. Lepowsky, Generalized Vertex Algebras and Relative
Vertex Operators, Progress in Math. {\bf Vol. 112}, Birkha\"user,
Boston, 1993.

\bibitem[DLM]{dlm}
C. Dong, H. Li  and G. Mason, Regularity of rational vertex 
operator algebras, {\em Adv. Math,} {\bf 132} (1997), 148-166.

\bibitem[DM]{dm}
C. Dong and G. Mason, On quantum Galois theory, {\em Duke Math. J.} 
{\bf 86} (1997), 305-321.

\bibitem[DN]{dn}
M. R. Douglas and N. A. Nekrasov, Noncommutative field theory,
arXiv: hep-th/0106048.

\bibitem[EK]{ek}
P. Etingof and D. Kazhdan, Quantization of Lie bialgebras, V,
arXiv: math.QA/9808121.


\bibitem[FHL]{fhl1}
I. Frenkel, Y.-Z. Huang and J. Lepowsky, On axiomatic approaches to
vertex operator algebras and modules,  Memoirs Amer. Math.
Soc. {\bf 104}, 1993; preprint, 1989.

\bibitem[FLM]{flm1}
I. Frenkel, J. Lepowsky and A. Meurman, {\it Vertex Operator Algebras
and the Monster}, Pure and Appl. Math. {\bf Vol. 134}, Academic Press,
Boston, 1988.

\bibitem[GL]{gl}
Y.-C. Gao and H.-S. Li, Generalized vertex algebras generated 
by parafermion-like operators, 
{\em J. Algebra} {\bf 240} (2001), 771-807.

\bibitem[K]{kac}
V. G. Kac, {\em Vertex Algebras for Beginners}, University Lecture Series, 
Vol. 10, 1996.

\bibitem[LL]{ll}
J. Lepowsky and H.-S. Li, Introduction to vertex operator algebras and 
their modules, Monograph, in preparation.

\bibitem[Li1]{li1}
H.-S. Li, Local systems of vertex operators, vertex superalgebras
and modules,
{\em J. Pure Appl. Alg.} {\bf 109} (1996), 143-195;  hep-th/9406185. 

\bibitem[Li2]{li2}
H.-S. Li, Local systems of twisted vertex operators, vertex superalgebras
and twisted modules, {\it Contemporary Math.} {\bf 193} (1996), 203-236.

\bibitem[Li3]{li3}
H.-S. Li, On certain higher dimensional analogues of vertex algebras,
arXiv:math.QA/0104225.

%\bibitem[Li4]{li4}
%H.-S. Li, Vertex (operator) algebra are ``algebras'' of vertex
%operators, in: Proceedings of the workshop on vertex operator 
%algebras, Toronto, Canada, October, 2000.

\bibitem[LZ]{lz}
B. Lian and G. Zuckerman, Commutative quantum operator algebras,
{\em J. Pure Appl. Alg.} {\bf 100} (1995), 117-139.

\bibitem[MN]{mn} 
A. Matsuo and K. Nagatomo, {\em Axioms for a Vertex Algebra and
the Locality of Quantum Fields}, MSJ Memoir, Vol. {\bf 4}, 
Mathematical Society of Japan, 1999.

\bibitem[Sn]{sn}
C. Snydal, Equivalence of Borcherds $G$-vertex algebras and axiomatic 
vertex algebras, math.QA/9904104.

\bibitem[Y]{y} 
G. Yamskulna, The relationship between skew group algebras and
orbifold theory, arXiv: math.QA/0107180.

\end{thebibliography}
\end{document}